\documentclass[11pt]{article}
\usepackage[margin=1.2in]{geometry}
\usepackage[utf8]{inputenc}
\usepackage[english]{babel}
\usepackage{amsmath,amsfonts,amssymb}
\usepackage{amsthm} 
\usepackage{hyperref}
\usepackage[capitalize]{cleveref}
\usepackage{graphicx}
\graphicspath{ {../} }
\usepackage{enumitem}
\usepackage{subfig}
\usepackage{capt-of}
\usepackage{mathdots}
\usepackage{mathrsfs}
\usepackage{authblk}
\makeatletter
\newcommand{\mylabel}[2]{#2\def\@currentlabel{#2}\label{#1}}
\makeatother
\usepackage{bm}
\usepackage{tikz,pgfplots}
\usetikzlibrary{fit,shapes,calc,backgrounds,patterns}
\usetikzlibrary{decorations.markings}
\usetikzlibrary {arrows.meta,bending} 
\usetikzlibrary{positioning,decorations.pathreplacing}
\usetikzlibrary{matrix,positioning}
\usepackage{bbding}
\usepackage{listings}
\usepackage{color}
\usepackage{csquotes}
\usepackage{mathtools}
\usepackage[most]{tcolorbox}
\usepackage{comment}
\usepackage{xcolor}
\usepackage{physics}
\usepackage{booktabs}
\usepackage{caption}
\usepackage[]{mdframed}
\usepackage{algorithm}
\usepackage{algpseudocode}
\newcounter{algsubstate}

\newtheorem*{theorem-non}{Theorem}

\theoremstyle{remark}

\theoremstyle{plain}
\theoremstyle{definition}

\newtheorem{example}{Example}

\newcommand{\M}{\pmb{\mathcal{M}}}
\newcommand{\vecsigma}{\pmb{\sigma}}
\newcommand{\R}{\mathbb{R}}
\newcommand{\mA}{\mathbf{A}}
\newcommand{\mB}{\mathbf{B}}
\newcommand{\mC}{\mathbf{P}}

\newcommand{\vecb}{\mathbf{b}}
\newcommand{\mS}{\mathbf{S}}

\newcommand{\mX}{\mathbf{X}}

\title{Data selection: at the interface of PDE-based inverse problem and randomized linear algebra}
\author{Kathrin Hellmuth\thanks{Department of Computing + Mathematical Sciences, Californian Institute of Technology, Pasadena, CA}, Ruhui Jin\thanks{Department of Applied Mathematics and Statistics, Johns Hopkins University, Baltimore, MD}, Qin Li\thanks{ Department of Mathematics and Wisconsin Institute for Discovery, University of Wisconsin-Madison,  Madison, WI}, Stephen J. Wright\thanks{Department of Computer Science, University of Wisconsin-Madison, Madison, WI}}

\date{}

\makeatletter
\patchcmd{\@maketitle}
  {\par\vskip 1.5em} 
  {\par\vskip 1.5em
   \begingroup
     
     \footnotetext{\textbf{Funding:}  QL, RJ and SW acknowledge support from the NSF through the grant NSF-DMS-2023239.  QL's work was also supported by NSF-DMS-2308440. KH acknowledges support from the Jet Propulsion Laboratory PDRDF 24AW0133.}
   \endgroup
  }{}{}
\makeatother

\begin{document}

\maketitle

\begin{abstract}
All inverse problems rely on data to recover unknown parameters, yet not all data are equally informative. This raises the central question of data selection. A distinctive challenge in PDE-based inverse problems is their inherently infinite-dimensional nature: both the parameter space and the design space are infinite, which greatly complicates the selection process. Somewhat unexpectedly, randomized numerical linear algebra (RNLA), originally developed in very different contexts, has provided powerful tools for addressing this challenge. These methods are inherently probabilistic, with guarantees typically stating that information is preserved with probability at least $1-p$ when using $N$ randomly selected, weighted samples. Here, the notion of “information” can take different mathematical forms depending on the setting. In this review, we survey the problem of data selection in PDE-based inverse problems, emphasize its unique infinite-dimensional aspects, and highlight how RNLA strategies have been adapted and applied in this context.
\end{abstract}

\section{Introduction}
Inverse problems based on partial differential equations (PDE) are relevant in such important and physically motivated applications as materials science, medical imaging, and geophysics. 
The fundamental task of the inverse problem is to reconstruct certain parameters (typically functions) in the PDEs that govern the dynamics of the system, using data gathered from the system. 
These parameters often encode critical information for downstream tasks. 
For example, recovering the optical properties of biological tissue may reveal abnormalities, recovering the refractive index of subsurface structures is essential in petroleum exploration, and recovering the band structure of a quantum system can help with material design.

PDE-based inverse problems have been studied extensively, from both analytical and algorithmic perspectives, across a wide range of PDE models (see, e.g., \cite{isakov, engl2000regularization, kirsch}). 
Most of this work assumes that data is given, and focuses on the quality and stability of reconstruction. 
In practice, reconstruction is often preceded by a critical step: {\em experimental design}.
Experimental design refers to the process of planning experiments so as to maximize the information content of the data collected for a given task. 
In PDE-based inverse problems, we usually cannot alter the governing equations themselves, but we do have control over experimental settings such as the source profile (which triggers the PDE dynamics) or the placement of detectors. 
This perspective is closely related to data selection: given the possibility of generating a large collection of data, the goal is to identify the most informative subset of experiments that yields the best possible reconstruction of the unknown parameters.

Classical experimental design is a rich and well-developed area of statistics, offering tools such as leverage scores, Fisher information, and optimality criteria including A-, D-, and E-optimality, which have been successfully applied across many scientific domains (see~\cite{F13, pukelsheim2006optimal}). 
In recent years, many of these techniques have been brought into the setting of PDE-based inverse problems, particularly within the framework of optimal design (see the reviews~\cite{A21, HJM24}). 
Yet, the direct transfer of experimental design principles to PDE-based problems is not straightforward; it introduces several distinctive challenges.
\begin{itemize}
\item {\em Structure induced by PDEs.} Classical design methods are largely generic and independent of specific models, but PDE-constrained inverse problems possess rich structural features, { commonly expressed as a tensorized structure,} that impose constraints while also offering opportunities. 
Effective experimental design must respect and, if possible, exploit these structures. 
\item {\em Efficiency over optimality.} Many design strategies seek an ``optimal" subset of experiments. 
In PDE-based settings, efficiency often takes precedence: Given the  {potentially} high  {computational costs associated to the PDE forward problem, the optimality requirement for data may be relaxed to some degree in favour of a more efficient data selection step. }  {This distinction lies at the heart of the disparity of quantitative design and qualitative design.}
\item {\em Cost of data acquisition.} Traditional experimental design often assumes that all candidate data are already available for evaluation. 
By contrast, in PDE-based inverse problems, each data point may be expensive to generate   {experimentally, whether in real world experiments or simulations, which}  
raises issues of sample complexity and calls for design strategies that account for acquisition costs. 
\end{itemize}
 These considerations indicate that classical methods for experimental design needed to be adapted in order to make them effective for PDE-based inverse problems.

The apparently unrelated research direction of {\em randomized numerical linear algebra} (RNLA) has developed rapidly in recent years~\cite{martinsson2020randomized,mahoney2011randomized}. 
Originally motivated by large-scale data problems in machine learning and internet applications such as recommendation systems, RNLA methods emphasize efficiency, trading some accuracy for scalability through randomization. 
The RNLA philosophy resonates with the goals of experimental design for PDE inverse problems, where efficiency is a more pressing concern than exact optimality.

One of the classical contributions of RNLA to experimental design is the fast computation of leverage scores in linear regression~\cite{DMMW12, mahoney2011randomized}, which serves as a cornerstone for randomized data selection. 
Over the years, these techniques have been adapted and extended, and are now used for experimental design in PDE-based inverse problems.

This review summarizes these recent advances at the interface of data selection for PDE-based inverse problems and RNLA. 
The field is still developing: Many algorithms are being actively refined, and nonlinear extensions remain largely open. 
Our aim is to provide a coherent overview of what has been achieved so far, while highlighting key challenges and promising directions for future research.

The remainder of the paper is structured as follows.
Section~\ref{sec:background_PDE} gives background on PDE-based inverse problems, highlighting structural features that are particularly relevant for experimental design. 
While reconstruction solvers are typically viewed as downstream tasks, the way they are formulated determines how much information can be extracted from data and thus directly guides data selection. 
Section~\ref{sec:preparation_tools} reviews two classical approaches: PDE-constrained optimization and Bayesian formulations. 
Both solvers rely critically on linearization; we also provide technical preparation in this context. 
Section~\ref{sec:rnla}, the main pillar of this review, surveys recent advances in RNLA methods and their applications to experimental design for PDE inverse problems. 
Section~\ref{sec:outlook} concludes with a discussion of nonlinear extensions and directions for future research.
\section{PDE based inverse problems}\label{sec:background_PDE}

Among inverse problems, those governed by PDEs form a distinct and important subclass, marked by unique structural and analytical challenges. 
While different PDEs give rise to diverse inversion properties, they also share fundamental traits that permit a unified treatment. 
A particularly notable feature of PDE-based inverse problems is their inherently infinite-dimensional nature, which contrasts with the finite-dimensional representations required for computation. 
It is not trivial to bridge this disparity.
Theoretical analysis is typically carried out in infinite-dimensional function spaces, whereas numerical algorithms must ultimately operate in finite dimensions. 
This bridge lies at the core of the present review and directly motivates our focus on data selection and experimental design.

In the subsections that follow, we introduce a generic PDE-based inverse problem in its infinite-dimensional form and describe the hierarchy of steps leading to the finite-dimensional numerical treatment. 
We discuss issues of well-posedness and the implications for experimental design, and we distinguish between quantitative and qualitative design goals. 
Finally, we conclude with representative examples of PDE-based inverse problems, grounding the abstract framework in concrete scientific applications.

\subsection{Infinite-dimensional nature}

Every PDE-based inverse problem begins with an underlying partial differential equation of the form:
\begin{equation}\label{eqn:PDE}
\mathcal{L}_{\vecsigma}[u] = S_{\theta_1}(x)\quad\Rightarrow\quad \mathrm{d}(\theta_1,\theta_2) = {M}_{\theta_2}[u]
\end{equation}
Here, $\mathcal{L}_{\vecsigma}$ denotes the PDE operator that models the physical or biological phenomenon under investigation, defined over a Euclidean domain with spatial variable $x \in \mathbb{R}^d$. The operator $\mathcal{L}_{\vecsigma}$ depends on a parameter $\vecsigma$, which encodes certain characteristics of the system. 
The input term $S_{\theta_1}(x)$ may represent a source term, boundary condition, or initial condition  {and is indexed by $\theta_1$}. 
The system is fully specified by fixing $\vecsigma$, and a given input $S_{\theta_1}(x)$ leads to a solution $u$ of the PDE.

The state $u$ is typically not observable in full. 
Instead, measurements are obtained through  $M_{\theta_2}$,  {  a functional that maps the PDE solution $u$ to a number that represents the reading of the solution through a detector indexed by $\theta_2$. }

 {Every pair 
$\theta = (\theta_1, \theta_2) \in \Theta$ 
specifies the experimental design configuration, with $\theta_1$ and $\theta_2$ respectively describing the design of source and detector, and gives rise to one} ground-truth datum $\mathrm{d}(\theta_1,\theta_2) = {M}_{\theta_2}[u]$.

Varying these configurations provides information about the unknown parameter $\vecsigma$, and the collection of all such configurations $\Theta$ is referred to as the design space.

This entire process defines what we refer to as the \textbf{Input-to-Output (ItO)} operator:
\begin{equation}\label{eqn:Lambda}
\Lambda_{\vecsigma}:\; (S_{\theta_1}, M_{\theta_2}) \mapsto \mathrm{d}(\theta_1, \theta_2)\,.
\end{equation}

Throughout this paper, we assume for all  { $\vecsigma$ in an admissible set $\Sigma$}, the PDE \eqref{eqn:PDE} is well-posed under a prescribed metric, and the solution is measurable using the observation operator ${M}_{\theta_2}$ for all values of $\theta_2$ under consideration. 
These assumptions guarantee that the ItO operator~\eqref{eqn:Lambda} is well-defined for each $\vecsigma$. Moreover, we equip the domain and codomain of $\Lambda_{\vecsigma}$ with metrics inherited from the input space $S_{\theta_1}(x)$ and the data space $\mathrm{d}(\theta_1,\cdot)$, respectively. Within this framework, the forward problem can be summarized as follows: for a fixed system parameter $\vecsigma$, the operator $\Lambda_{\vecsigma}$ maps each input $S_{\theta_1}(x)$ to a corresponding data function $\mathrm{d}(\theta_1,\cdot)$.

The inverse problem is to recover the parameter $\vecsigma$  {living in a parameter space $\Sigma$} from knowledge of the operator $\Lambda_{\vecsigma}$, or equivalently, from the full collection of data values $\mathrm{d}(\theta_1,\theta_2)$ over $\theta = (\theta_1,\theta_2) \in \Theta$. To formalize this, we define the data map:
\begin{equation}\label{eqn:data}
\M(\theta; \vecsigma): \Theta \otimes \Sigma \to \mathbb{R},\;\quad \text{with} \quad \M(\theta, \vecsigma) = \Lambda_{\vecsigma}[S_{\theta_1}, M_{\theta_2}] = \mathrm{d}(\theta_1, \theta_2)\,.
\end{equation}
The task of the inverse problem is then to infer $\vecsigma$ from complete knowledge of $\M(\cdot\;;\vecsigma)$.

A defining characteristic of this class of problems is their infinite-dimensional nature. The unknown parameter $\vecsigma$ is typically a function and thus infinite-dimensional, while the data function $\mathrm{d}$ is indexed continuously over the design space $\Theta$. Consequently, both the parameter space $\Sigma$ and the design space $\Theta$ 
reside in infinite-dimensional settings, distinguishing these problems from standard finite-dimensional inverse problems.

\subsection{Finite-dimensional computation}
While the natural formulation of PDE-based inverse problems is often infinite-dimensional, we need to reduce to finite dimensions for experimental and numerical reasons.
This reduction arises from two primary sources:
\begin{itemize}
    \item \textbf{Limited data:} Only finitely many experiments can be carried out in practice, which restricts the design space to a finite subset $\Theta_c \subset \Theta$. 
    If $c$ experiments are conducted, and each experiment requires a specific value for $\theta=(\theta_1,\theta_2)$, then we are taking $|\Theta_c| = c$.    
    \item \textbf{Reduced parameterization:} 
     $\vecsigma$ is typically represented with a finite number of degrees of freedom, yielding a finite-dimensional approximation $\Sigma_n$ of dimension $n$. Any $\vecsigma$ living in this space $\Sigma_n$ is then automatically a finite-dimensional object.
   \end{itemize}
Under these approximations, the data map~\eqref{eqn:data} becomes
\begin{equation}\label{eqn:data_finite}
    \M(\theta; \vecsigma):\; \Theta_c \otimes \Sigma_n \to \mathbb{R}\,.
\end{equation}

It is tempting to view the finite-dimensional formulation as simply an approximation of the infinite-dimensional one. Indeed, as $n,c\to\infty$, the two settings coincide. However, theoretical results established in the infinite-dimensional framework do not in general translate directly to the finite $n,c$ setting. Bridging these regimes requires a hierarchy of intermediate analyses, which we summarize as follows.
\begin{itemize}
    \item \textbf{Well-posedness at the PDE level:} In the infinite-dimensional setting, well-posedness concerns uniqueness and stability of the inverse map. 
    We seek an estimate of the form  
\begin{equation}\label{eqn:well_posed}
    \|\vecsigma_1 - \vecsigma_2\| \leq w\!\left(\|\M(\cdot, \vecsigma_1) - \M(\cdot, \vecsigma_2)\|\right)\,, \quad  {\vecsigma_1,\vecsigma_2} \in \Sigma \cap B_R(\vecsigma_\ast)\,,
\end{equation}  
where  { the modulus of continuity $w$ of $\mathcal{M}$} quantifies the stability of the inverse map and $B_R(\vecsigma_\ast)$ denotes the neighborhood of radius $R$ around the ground truth $\vecsigma_\ast$.  Here $R$ may be infinite, or may encode problem-specific constraints such as positivity. 
 {Naturally, the properties of $w$  depend strongly on the choice of topology on the space, and thus the norm $\|\cdot\|$ used in~\eqref{eqn:well_posed}.} The following two properties capture the essence of well-posedness.
    \begin{itemize}
        \item \textbf{Uniqueness:} \begin{equation}\label{eqn:uniqueness}
    w(0) = 0\,,
\end{equation}  
implying that complete knowledge of $\M(\cdot,\vecsigma)$ uniquely determines $\vecsigma$. 
Uniqueness has been established for many PDE-based inverse problems, the most celebrated being the Calder\'on problem~\cite{calderon},  {see \ref{itm: E2} below}, with pioneering results in~\cite{uhlmann_annals,kohn1984determining} and subsequent extensions in~\cite{novikov1988multidimensional,colton1992uniqueness,isakov1993uniqueness,nachman}.

        \item \textbf{Stability:} We seek constants $C > 0$ and $\beta > 0$ such that  
\begin{equation}\label{eqn:stability}
    w(t) \leq C t^\beta\,,
\end{equation}  
ensuring that small perturbations in data lead to controlled (e.g., H\"older/Lipschitz continuous) perturbations in the reconstruction. 
It is also possible for the upper bound to take on a logarithmic or exponential form. 
Foundational results in this direction include inverse conductivity problems~\cite{alessandrini1988stable} and inverse scattering problems~\cite{Hahner_stability,stefanov1990stability,bellout1992stability}. 
    \end{itemize}
These results are highly problem-specific, relying on PDE-dependent techniques. For instance, complex geometrical optics (CGO) solutions are central in Schr\"odinger-type and elliptic problems~\cite{uhlmann_salo}, Carleman estimates play a key role in transport and hyperbolic equations~\cite{KLIBANOV2008352,Klibanov_1992,Yamamoto}, and singular decomposition or averaging-lemma techniques are employed for kinetic equations~\cite{choulli_stefanov,Li_2020}.  
    
    \item \textbf{Well-posedness in the finite-dimensional setting:} After discretization, we ask analogous questions of uniqueness and stability with $\Sigma$ replaced by $\Sigma_n$ and $\Theta$ by $\Theta_c$. The relevant estimate becomes  
\begin{equation}\label{eqn:well_posed_finite}
    \|\vecsigma_1 - \vecsigma_2\| \leq w_{n,c}\!\left(\|\M(\cdot|_{\Theta_c}, \vecsigma_1) - \M(\cdot|_{\Theta_c}, \vecsigma_2)\|\right)\,, \quad \sigma \in \Sigma_n \cap B_R(\vecsigma_\ast)\,.
\end{equation}  
That is:
\begin{itemize}
\item If $w_{n,c}(0) = 0$, then $\M(\cdot|_{\Theta_c}, \vecsigma)$ uniquely determines $\vecsigma$ within $\Sigma_n$.  
\item If an inequality of the form~\eqref{eqn:stability} holds for $w_{n,c}$, then small data perturbations lead to bounded reconstruction error.  
\end{itemize}
It is worth noting that well-posedness in the infinite-dimensional setting does not directly translate to the finite-dimensional case. The properties of $w$ in~\eqref{eqn:well_posed} are only suggestive; they do not guarantee corresponding properties for $w_{n,c}$ in~\eqref{eqn:well_posed_finite}.   {In particular, the properties of $\omega_{n,c}$ crucially depend on the choice of $\Sigma_n$ and $\Theta_c$.}
    \item \textbf{Numerical implementation:} Ultimately, we need to implement a computational scheme for the reconstruction. This task requires algorithmic pipelines that recover $\vecsigma$ from observed data. Two broad approaches dominate: (i) PDE-constrained optimization methods~\cite{hinze2008optimization,haber2000optimization}, and (ii) Bayesian inference frameworks~\cite{kaipio2005statistical,stuart2010inverse}. 
    When prior information is available, regularization methods are widely used to stabilize reconstruction, including classical Tikhonov regularization~\cite{1571980074135178112,ito2014inverse}, total variation (TV)~\cite{rudin1992nonlinear}, and $\ell_1$-based sparsity methods~\cite{daubechies2004iterative,jin2012sparsity}.  
\end{itemize}
Throughout this paper, we assume that PDE-level well-posedness results~\eqref{eqn:well_posed}--\eqref{eqn:stability} are available. Our focus is on transferring these theoretical guarantees to the finite-dimensional setting, establishing such finite-dimensional counterparts as~\eqref{eqn:well_posed_finite}.

\subsection{Well-posedness for the finite-dimensional problem and design framework}
\label{sec:infinite_to_finite}
\begin{figure}[hbt]
    \centering
    \includegraphics[width=0.5\linewidth]{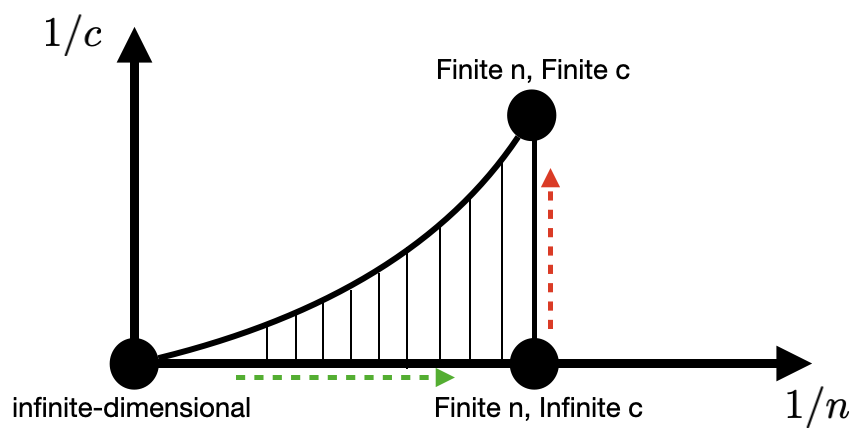}
    \caption{Conceptual pathway from infinite- to finite-dimensional formulations. The origin $(n,c) = (\infty,\infty)$ represents the infinite-dimensional setting, where well-posedness is established. The green arrow reduces the parameter space to finite dimension while keeping data abundant, leading to the semi-infinite case. The red arrow then reduces the data to finitely many experiments, posing the key challenges of identifiability and data selection. }
    \label{fig:from_infinite_to_finite}
\end{figure}

It is tempting to transfer the infinite-dimensional well-posedness estimate~\eqref{eqn:well_posed} directly to its finite-dimensional counterpart~\eqref{eqn:well_posed_finite}. The passage, however, is not straightforward. 
A conceptual illustration appears in Figure~\ref{fig:from_infinite_to_finite}.

At the origin $(n,c) = (\infty,\infty)$, the problem is posed in the infinite-dimensional setting, where well-posedness is assumed to hold. The most natural route to the fully finite-dimensional case proceeds first along the green arrow, followed by the red arrow.

\medskip
\noindent\textbf{The green arrow.} This step reduces the problem to a semi-infinite setting: the unknown parameter is finite-dimensional, while data remains abundant. This reduction is essentially free. Once well-posedness~\eqref{eqn:well_posed} is established in the infinite-dimensional formulation, uniqueness readily follows in the semi-infinite case, with
\begin{equation}\label{eqn:well_posed_semifinite}
|\vecsigma_1 - \vecsigma_2| \leq w_n\left(|\M(\cdot, \vecsigma_1) - \M(\cdot, \vecsigma_2)|\right), \quad \text{for } \sigma \in \Sigma_n \cap B_{R}(\vecsigma_\ast),
\end{equation}
for some function $w_n$. 
Since the amount of data used in the same, but the to-be-reconstructed parameter $\vecsigma$ now lives in a smaller  {set}: $\Sigma_n\subset\Sigma$, the stability function $w_n$ inherits $w_n(0)=0$ from $w$. Moreover, because $\Sigma_n$ is finite dimensional and the analysis is restricted to a neighborhood of the ground truth, the domain under consideration is essentially compact, and the mapping is automatically Lipschitz~\cite{BOURGEOIS2013187}. 
 {We should stress that $\Sigma_n$ is  not necessarily a linear subspace. Indeed, setting $\vecsigma$ as a Gaussian mixture is a common example, in which the induced $\Sigma_n$ is a nonlinear manifold within $\Sigma$, formed by a set of parameters indicating means and variances of the Gaussians.}

Well-posedness in this semi-infinite regime has been studied extensively. For example, under the a priori assumption of piecewise-constant parameters, Alessandrini and Vessella~\cite{alessandrini2005lipschitz} established Lipschitz stability for the Calder\'on problem~\cite{calderon}. Subsequent work has extended these results to more general settings~\cite{alberti2019calderon,beretta2011lipschitz,gaburro2015lipschitz,ruland2018lipschitz,alberti2019calderon,alessandrini2017lipschitz,harrach2019uniqueness,Bacchelli_2006}, as well as to Schr\"odinger and Helmholtz equations~\cite{doi:10.1137/120869201,beretta2016inverse,alessandrini_deHoop_lipschitz}. However, the Lipschitz stability constants may grow potentially exponentially with the parameter dimension $n$~\cite{rondi2006remark,beretta2013lipschitz,li2021lipschitz}.

\medskip
\noindent\textbf{The red arrow.} In this (far more challenging) step, the available data is reduced from unlimited ($c=\infty$) to a finite $c$, while one seeks to maintain a comparable quality  of reconstruction. This reduction raises two fundamental questions.
\begin{enumerate}
\item[\mylabel{itm: Q1}{{\bf Q1.}}] What is the minimal number of experiments $c$ required to guarantee identifiability?
\item[\mylabel{itm: Q2}{{\bf Q2.}}] How should the subset $\Theta_c\subset\Theta$ be chosen?
\end{enumerate}
These questions lie at the heart of {data selection}, which is essential for bridging continuous formulations with practical, finite data. 
Addressing them is central to PDE-based inverse problems and is the primary focus of this review.

Finally, one may also define a particular structure for $\Theta_c$ and examine only the scaling of $c$ in terms of $n$. 
This approach was taken in a recent paper~\cite{arma_alberti}, where $\Theta_c$ is assumed to be a $c$-dimensional projection with a generic assumption.

\subsection{Quantitative vs. Qualitative}
The task of selecting data from a potentially infinite pool of experimental configurations is classical, and many approaches have been developed to address it. 
Broadly speaking, these approaches can be grouped into two categories.
\begin{itemize}
    \item[--]\underline{Quantitative design.} This perspective is rooted in the Bayesian framework of \emph{optimal experimental design}. Here, ``design" refers to selecting, from all possible configurations, the experiments that most effectively extract information about the unknown parameter $\vecsigma$. A quantitative design requires specifying a criterion $\Phi$ that measures the information content of a chosen subset $\Theta_c\subset\Theta$, then choosing the design that optimizes $\Phi$, that is 
    \begin{equation}\label{eq:OEDobjective}
\Theta_c^* = \arg \min_{\Theta_c\subset \Theta} \Phi(\Theta_c)\,.
\end{equation}
Over the past decades, a wide variety of criteria  {$\Phi$} and associated computational strategies have been developed. 
Among them, criteria based on minimizing uncertainty are particularly prominent: one models perturbations in the data and studies how they propagate into the reconstruction of \(\vecsigma\), then selects experiments that minimize the resulting uncertainty. 
See Section~\ref{sec:bayes} for further details. 
For broader treatments of experimental design we refer to classical solvers~\cite{W70,A73,MM70,W72} and to the summaries in~\cite{F13,pukelsheim2006optimal, drovandi, rainforth2024modern}. 
Design used in PDE-based inverse problems context is reviewed in~\cite{A21, HJM24,Chen_Santosa_2024}.

\item \underline{Qualitative design.} 
More recent approaches move away from strict optimality criteria. 
Instead, the goal is to secure reconstructions of reasonable quality \cite{walter1990qualitative, walter1996identifiability}. 
This relaxed perspective broadens the toolbox: techniques from randomized numerical linear algebra, gradient flows, and related areas that are traditionally viewed as separate from optimal design can now be applied effectively.
\end{itemize}

While quantitative design has long been the dominant paradigm, qualitative design has recently gained traction, particularly in connection with randomized numerical linear algebra (RNLA). 
As we discuss in Section~\ref{sec:rnla}, RNLA  prioritizes efficiency while tolerating modest losses in precision. 
This philosophy aligns naturally with qualitative design, where ``good enough" data selection suffices for practical purposes. 
For this reason, we are chiefly interested in qualitative approaches in this review.

\subsection{Some representative PDE-based inverse problems}\label{sec:representatives}

Several inverse problems arising from physical and biological applications have become canonical in the study of PDE-based inverse problems. 
We briefly summarize  {four} such examples here.
They are the inverse elliptic equation (also known as electrical impedance tomography, or the Calder\'on problem), the inverse Schr\"odinger equation, the inverse transport equation (as in optical tomography), and the inverse Helmholtz equation (arising in acoustic imaging).
 {As the terminology and methodology also apply to ODE-based inverse problems, we include an ODE example as Example E1.}

\begin{enumerate}
    \item [\mylabel{itm: E1}{{\bf E1.}}] \textbf{Inverse Lorenz 63 system.} 
The Lorenz 63 system \cite{DeterministicNonperiodicFlow}, originally introduced as a simplified model of atmospheric convection, describes the temporal evolution of three state variables  {in terms of the ODEs}:
\begin{equation}
    \begin{cases}
        \dot{x} = \sigma(y - x)\,, \\
        \dot{y} = x(\rho - z) - y\,, \\
        \dot{z} = xy - \beta z\,.
    \end{cases}
\end{equation}
Here, $(x,y,z)$ correspond to convective flow rate, horizontal temperature difference, and vertical deviation of the temperature profile, while the parameter vector $\vecsigma=(\sigma,\rho,\beta)$ encodes the Prandtl number, the Rayleigh number, and a geometric factor. Owing to its chaotic dynamics and bifurcation structure, the Lorenz system is a classical model in nonlinear dynamics. 

In the inverse setting, given partial or full time-series observations of $(x,y,z)$, the task is to reconstruct the parameter vector $\vecsigma$. The associated ItO map is:
\[
\Lambda_{\vecsigma} : (x, y, z)(t=0) = \theta_1 \mapsto (x, y, z)(t = \theta_2),\,.
\]
 {where $\theta_1$ and $\theta_2$ describe initial condition, and measurement time.}
Parameter identification for the Lorenz system has been widely studied; see, for instance, \cite{10.1063/1.1635095, KUNZE2007124, shatalov}.

\item [\mylabel{itm: E2}{{\bf E2.}}] \textbf{Inverse elliptic equation.} 
Elliptic equations provide the foundation of steady-state modeling. A typical formulation is:
\begin{equation}\label{eqn:elliptic}
\nabla \cdot (\vecsigma(x) \nabla u) = 0\,.
\end{equation}
The forward problem maps a prescribed Dirichlet boundary condition 
\begin{equation}\tag{\theequation i}\label{eqn:elliptic:DirichletBoundary}
    u|_{\partial \Omega} = \phi
\end{equation} 
to the solution $u$. A canonical inverse problem is electrical impedance tomography (EIT) \cite{EIT_review}, a noninvasive medical imaging technique in which one seeks to infer the conductivity distribution $\vecsigma$ of tissue from boundary voltage–current measurements. In this case, the measured Neumann data $\vecsigma\partial_n u$ defines the Dirichlet-to-Neumann (DtN) map {, that we use as our ItO map in the data selection setting}:
\[
\Lambda_{\vecsigma} : \phi = \phi_{\theta_1} \mapsto \vecsigma \, \partial_n u(\theta_2)\,
\]
 {with $\theta_1$ parameterizing the Dirichlet boundary condition and  $\theta_2$ denoting the  measurement location of the Neumann data.}
Physically, this corresponds to the  voltage-to-current relation.  {There are a few variations to this problem. One important variation is to deploy the current-to-voltage data for the reconstruction, a map that is mathematically equivalent to the Neumann-to-Dirichlet map:
\[
(\Lambda_{\vecsigma})^{-1} : \psi_{\theta_1} \mapsto u(\theta_2)\,.
\]
This means~\eqref{eqn:elliptic} is equipped with Neumann boundary data $\vecsigma \partial_n u = \psi_{\theta_1}$ that describe the current on $\partial \Omega$, indexed by $\theta_1$, and Dirichlet data are taken on boundary points $\theta_1\in \partial \Omega$ that represents the voltage measurements.}  {Another variation is usually termed the Darcy flow problem~\cite{HV19,GHLS20}. Mathematically this is to add a source term to the elliptic equation $\nabla \cdot (\vecsigma(x) \nabla u) = S(x)$, and take measurement directly on $u$ in the interior. By varying the source $S(x)$ and the detector location for $u$, one infers $\vecsigma$. The problem is mathematically simpler, but requires intrusive measurements.}

The EIT problem has a rich history; see surveys \cite{borcea2002electrical, uhlmann30, uhlmann_salo, Adler2020}. Foundational works established uniqueness and stability results \cite{uhlmann_annals, kohn1984determining, alessandrini1988stable, nachman, kari}.

        \item [\mylabel{itm: E3}{{\bf E3.}}]\textbf{Inverse Schr\"odinger equation} The Schr\"odinger equation is a fundamental model for describing wave fields in quantum systems. Its steady-state form is a  {linear} elliptic equation: \begin{equation}\label{eqn:schroedinger}
 \Delta_x u  {-} \vecsigma(x)u = 0\,, \quad \text{with} \quad u|_{\partial \Omega} = \phi\,.
 \end{equation}
 where $u$ is the wave field and $\vecsigma$ the potential. The inverse problem seeks to reconstruct $\vecsigma$ from boundary measurements of $u$. Applications include material characterization and design for achieving prescribed properties~\cite{PhysRevResearch.7.013296,molecule_design_science,inverse_band}. The corresponding boundary map,  {that we use to construct ItO map,} is:
 \[
 \Lambda_{\vecsigma}: {\phi_{\theta_1}\mapsto
  u (\theta_2)} \,,
 \]
  {with $\theta_{1}$ carrying  a similar meanings as in \ref{itm: E2}, and $\theta_2$ denoting a measurement location within the domain.}
 From a mathematical perspective, the inverse Schr\"odinger problem can be related to EIT: Under smoothness assumptions on the medium, the Schr\"odinger potential ($\vecsigma$ in~\eqref{eqn:schroedinger}) can be written as $\tfrac{\Delta \sqrt{\vecsigma}}{\sqrt{\vecsigma}}$ from~\eqref{eqn:elliptic}, linking the two formulations~\cite{nachman, uhlmann_salo}.

    \item [\mylabel{itm: E4}{{\bf E4.}}]\textbf{Inverse transport equation.} Transport equations are widely used in optical and radar imaging. The underlying principle is to infer the optical properties of a medium from measurements of incoming and outgoing photon fluxes~\cite{OT_human_brain,FERRARI2012921,Bluestone:01,Mourant:98}. These material properties are encoded in the coefficient $\vecsigma$. The governing PDE is
    \[
    v \cdot \nabla_x f = \vecsigma(x) \, \mathcal{L}[f]\,, \quad \text{with} \quad f|_{\Gamma_-} = \phi(x, v)\,,
    \]
    where $f(x, v)$ denotes the photon distribution at spatial location $x$ and velocity $v$, and $\mathcal{L}$ is a prescribed scattering operator. The domain is $x \in \Omega$, with $\Gamma_\pm = \{ (x, v) : x \in \partial \Omega, \pm v \cdot n_x > 0 \}$ representing inflow and outflow boundaries. The associated inverse problem seeks to reconstruct $\vecsigma$ from boundary measurements, which defines the ItO map:
    \[
    \Lambda_{\vecsigma} : \phi_{\theta_1} \mapsto f|_{\Gamma_+}(\theta_2)\,.
    \]
     {In many applications, laser beams used to inject photons are very focused, and $\theta_1$ codes the laser beam location and direction on the boundary. Similarly, $\theta_2$ codes the location of the detector and the directional information the detector can receive.}
    Inverse transport problems are fundamental in optical tomography and have been studied extensively from both analytical and computational perspectives \cite{larsen_inverse_transport, choulli_stefanov, Wang1999, Bal_2009, Arridge_2009, ren2010recent}.

    \item [\mylabel{itm: E5}{{\bf E5.}}]\textbf{Inverse Helmholtz equation.} The Helmholtz equation models time-harmonic wave propagation and forms the foundation of acoustic imaging, with applications in geophysics~\cite{bergen2019machine,Matt_Li_MMS,Demanet} and medical imaging~\cite{rosenthal2013acoustic}. The PDE takes the form
    \[
    \Delta_x u + \omega^2 \vecsigma(x) u = 0\,,
    \]
    where $\omega$ is the prescribed frequency.  The incoming wave is specified as a plane wave $u^{\text{in}}_{\theta_1} = e^{i \omega \theta_1 \cdot x}$  {propagating in the direction of $\theta_1$}.  {Then the total field $u$ decomposes into $u = u^{\text{in}}+u^{\text{sc}}$ into the incoming wave and scattered wave $u^{\text{sc}}$ emerging from the interaction with the scatterer described by $\vecsigma$}. Measurements are collected on the boundary of a large domain (typically modeled as a ball of radius $R$, denoted $B_R$). The inverse problem consists of reconstructing the coefficient $\vecsigma$ using the ItO map:
    \[
    \Lambda_{\vecsigma} : u^{\text{in}}_{\theta_1} \mapsto u|_{\partial B_R}(\theta_2)\,,
    \]
     {where $\theta_2$ codes the location of the detector on the ring of $\partial B_R$.}  {The domain size $R$ determines the inversion regime. It is common to use far-field data, meaning $R\gg 1$. }
    The inverse Helmholtz scattering problem has been the subject of extensive research from mathematical, statistical, and numerical perspectives. 
    Works that summarize the progress and highlight key advances in theory and algorithm design include the books~\cite{kirsch, fouque2007wave, colton2012inverse, ammari, borcea2014imaging,isakov} and landmark works~\cite{Tarantola_seismic, borcea2002imaging, ito2012direct, bao2015inverse, engquist2011sweeping, taus2020sweeps,borcea2014imaging}.
\end{enumerate}

\section{Technical preparation and classical tools}\label{sec:preparation_tools}
Experimental design is an upstream task that precedes the implementation and execution of inverse algorithms. The value of the collected data ultimately depends on how effectively it can be exploited by the chosen solver. Two major algorithmic pipelines have been particularly influential: \emph{PDE-constrained optimization} and the \emph{Bayesian formulation}, both of which we summarize in the following subsections. A recurring theme in both frameworks is the central role of \emph{linearization} that requires gradient computation. We therefore devote Section~\ref{sec:linearization} to surveying the linearization techniques.

\subsection{PDE-constrained optimization}\label{sec:PDEOptimization}
When data is abundant, PDE-constrained optimization is one of the most widely used methods for inferring unknown parameters in PDEs. The formulation is straightforward: Subject to the PDE being satisfied, we seek the parameter $\vecsigma$ that produces the best match to the measurement data:
\begin{equation}\label{eqn:PDE_constrained_opt}
    \min_{\vecsigma}\sum_{\theta_1,\theta_2}\frac{1}{2}|\mathrm{d}(\theta)-M_{\theta_2}[u_{\theta_1}]|^2
    \,,\quad\text{s.t.} \quad \mathcal{L}_{\vecsigma}[u_{\theta_1}]=S_{\theta_1}\,.
\end{equation}
If the data $\mathrm{d}$ is noiseless (as when it is generated by the PDE solution using a ground-truth parameter $\vecsigma^*$) and the problem is well-posed, then the optimizer is unique and yields a zero objective value. Since the PDE itself is well-posed, $u_{\theta_1}$ is uniquely determined; we therefore write it as $u(\vecsigma;\theta_1)$ to emphasize its dependence on $\vecsigma$. 
The optimization problem can then be recast  {as follows:}
\[
\min_{\vecsigma}L(\vecsigma)=\frac{1}{2}\sum_{\theta_1,\theta_2}|\mathrm{d}(\theta)-M_{\theta_2}[u(\vecsigma;\theta_1)]|^2=\frac{1}{2}\sum_{\theta_1,\theta_2}|\mathrm{d}(\theta)-\mathrm{pde}(\vecsigma,\theta))|^2\,,
\]
where we adopt the shorthand $\mathrm{pde}(\vecsigma,\theta) = M_{\theta_2}[u(\vecsigma;\theta_1)]$.

In practice, gradient-based solvers are most commonly used to minimize $L(\vecsigma)$, updating $\vecsigma$ by making use of the gradient $\frac{d}{d\vecsigma}L(\vecsigma)$. 
The conditioning of the problem, which governs both stability and convergence properties,
depends critically on the Hessian $\mathrm{Hess}_{\vecsigma}L$. 
Calculation of the gradient and Hessian  can be summarized as follows.
\begin{enumerate}
    \item \textbf{Gradient} $\frac{d}{d\vecsigma}L$. By the chain rule, we have
    \[
    \frac{d}{d\vecsigma}L = \sum_\theta \left(\mathrm{pde}(\vecsigma,\theta)-\mathrm{d}(\theta)\right)\frac{d}{d\vecsigma}\mathrm{pde}(\vecsigma,\theta)\,.
    \]
    Since the dependence on $\vecsigma$ appears only through $\mathrm{pde}(\vecsigma,\theta)$, computing the gradient amounts to evaluating $\tfrac{d}{d\vecsigma}\mathrm{pde}$.
    \item \textbf{Hessian $\mathrm{Hess}_{\vecsigma}L$.}  
Differentiating once more gives  
\begin{equation}\label{eqn:hessian_L}
\mathrm{Hess}_{\vecsigma}L \;=\; 
\underbrace{\sum_{\theta}\frac{d}{d\vecsigma}\mathrm{pde}(\vecsigma,\theta)\otimes \frac{d}{d\vecsigma}\mathrm{pde}(\vecsigma,\theta)}_{\text{Fisher Information Matrix}} 
+ \sum_{\theta}\bigl(\mathrm{pde}(\vecsigma,\theta) - \mathrm{d}(\theta)\bigr)\,\mathrm{Hess}_{\vecsigma}\mathrm{pde}(\vecsigma,\theta)\,.
\end{equation}  
The Hessian thus consists of two components. The first is the outer product of the pde gradients, usually referred to as the Fisher Information Matrix (FIM) or the Gauss-Newton matrix. Near the ground truth, where $\mathrm{pde}(\vecsigma,\theta)\approx \mathrm{d}(\theta)$, the second term vanishes, leaving the Hessian dominated by the FIM. Note that when $\vecsigma \in \Sigma$, a function space, the Hessian is an operator, whereas for $\vecsigma \in \Sigma_N$ it reduces to a finite-dimensional matrix in $\mathbb{R}^{N\times N}$.  
\end{enumerate}
Both the gradient and the Fisher Information Matrix (FIM) will reappear as key quantities in experimental design, since they provide the sensitivity information that guides the selection of data so as to maximize its impact on parameter inference.

\subsection{Bayesian formulation}\label{sec:bayes}
Bayes’ formula provides a widely used alternative to PDE-constrained optimization, particularly when we need to quantify uncertainty of reconstruction. 
Instead of seeking a single best-fit parameter, the Bayesian perspective treats a range of parameters as plausible, each weighted by its posterior probability.

In this framework, prior knowledge is blended with observed data. For PDE-based inverse problems, we model the data as
\[
\mathrm{d}(\theta)=\mathrm{pde}(\vecsigma,\theta)+\eta(\theta)\,,
\]
where $\eta$ represents measurement error. The posterior distribution of $\vecsigma$ --- which is also the conditional distribution of $\vecsigma$ given the data $\mathrm{d}$ --- takes the form
\[
p^\mathrm{d}(\vecsigma) = p(\vecsigma|\mathrm{d})\propto p(\mathrm{d}|\vecsigma)p(\vecsigma) = p_\eta\left(\mathrm{d}(:)-\mathrm{pde}(\vecsigma,:)\right)p(\vecsigma)\,,
\]
where $p(\vecsigma)$ is the  {\em marginal} or {\em prior} distribution distribution and $p_\eta$ is the noise distribution. 

As in optimization, the posterior distribution has a few key quantities, one of which is the variance, which measures the uncertainty of the reconstruction. 
In general, the variance is difficult to compute for nonlinear PDE operators. 
However, in the linearized setting,  {the problem falls in the Gaussian regime, and the posterior distribution is well understood. In this regime,} the following assumptions are standard.
\begin{itemize}
    \item The PDE operator is linear, e.g. $\mathrm{pde}(\vecsigma,\theta)=\mA_{\theta,:}\vecsigma$, or locally linearizable around $\vecsigma^*$ as $\mathrm{pde}(\vecsigma,\theta)\approx \mA_{\theta,:}\left(\vecsigma-\vecsigma^*\right)+\mathrm{pde}(\vecsigma^*,\theta)$, with data adjusted accordingly.
\item Noise is Gaussian: $\eta\sim\mathcal{N}(0,\Gamma)$;
\item The prior is Gaussian: $\vecsigma\sim\mathcal{N}(\vecsigma^\dagger,\Gamma_0)$.
\end{itemize}
Under these assumptions, the posterior remains Gaussian, { and its covariance matrix is 
$$
    \hat \Gamma  =(\mA^\top \Gamma^{-1}\mA+\Gamma_0^{-1})^{-1}\,.$$
}
 {In analogy to}  the optimization framework, the matrix $\mA$, which captures the linear component of $\mathrm{pde}(\vecsigma,\theta)$, plays a central role.

\subsection{Bayesian Optimal Experimental Design}
 {
In classical literature, one measures information of data $\mathrm{d}(\theta)$ w.r.t. the parameter $\vecsigma$ by checking its contribution in the reconstruction uncertainty. Smaller uncertainty is preferred for a more accurate and confident reconstruction, and data should be selected accordingly. Typically the quantification of such uncertainty is expressed by the Bayes formula, making it one of the leading approaches used in quantitative data selection, and is paired with~\eqref{eq:OEDobjective}.

Use linear setting described above as an example, the covariance matrix of the posterior distribution offers a natural candidate for quantifying uncertainties, and a subset of data $d|_{\Theta_c}$ contribute the uncertainties coded in $\hat \Gamma|_{\Theta_c}$. As a consequence, to conduct data selection is to find $\Theta_c$ so that $\hat \Gamma|_{\Theta_c}$ is small. Since $\hat \Gamma|_{\Theta_c}$ is a matrix, some scalarization is needed, and popular choices are~\cite{kieferADoptimality1959,mitchell2000algorithm,alexanderian2018efficient,bandara2009optimal,park2018optimal,attia2022optimal,aarset2025global, pukelsheim1993optimal}:
\begin{itemize}
    \item[--] E-optimal design: $\Phi_E {(\Theta_c)} =  \lambda_{\max}(\hat \Gamma|_{\Theta_c})$.
    \item[--] D-optimal design: $\Phi_D(\Theta_c) = \det(\hat \Gamma|_{\Theta_c})$.
    \item[--] A-optimal design: $\Phi_A(\Theta_c) = \trace(\hat \Gamma|_{\Theta_c})$.
\end{itemize}
Accordingly, the Bayesian optimal experimental design refers to minimize these uncertainty measures:
\begin{equation}\label{eqn:optimal_design}
\min_{\{\Theta_c\subset \Theta\}} \;\Phi_{E,D,A}({\hat \Gamma|_{\Theta_c}})\,.
\end{equation}
This is a combinatorial problem due to the binary nature of the formulation. In practice, one often resorts to weights assignment as a relaxed version of the formulation. When the data point $\mathrm{d}(\theta)$ is weighted by $w_\theta\in \mathbb R$,
the posterior distribution becomes
\begin{equation}\label{eq:Bayesian inversion}
p^\mathrm{d, W}(\vecsigma)\sim\mathcal{N}(\hat \vecsigma_W, \hat \Gamma_W)\, \quad \text{with}\quad \hat \Gamma_W  =(\mA^\top W^{1/2}\Gamma^{-1}W^{1/2}\mA+\Gamma_0^{-1})^{-1}\,,
\end{equation}
Here $W=\text{diag}\{w_{\theta}\}$. This reformulation now allows one to translates the afore mentioned design criteria to be functions of $W$, and the optimization becomes:
\begin{equation*}\label{eqn:optimal_design_weight}
\min_{W=\{w_\theta\}} \;\Phi_{E,D,A}(W)\,,
\end{equation*}
The optimizer assigns higher weights to data most effective at reducing uncertainty. In many contexts, only a finite budget $c$ of experiments can be realized, and it is a convention to collects the experimental configurations associated with the $c$ hightest weights $w_\theta$, see~\cite{bouhtou2010submodularity}. To ensure a clear selection of a finite number of experiments, it is also common to add a regularizer that promotes sparsity of the weights $\{w_\theta\}$ during the optimization~\cite{koval2020optimal,haber2008numerical,APSG14}.

It is important to note that in the semi-infinite case, where $\vecsigma \in \Sigma_N$ and $\hat\Gamma_W \in \mathbb{R}^{N\times N}$, the sequence ${w_\theta}$ may still be infinite in size, meaning the optimization is still conducted over the entire design function space $\Theta$.  

The methodology can also be extended to nonlinear non-Gaussian regime. New challenges appear, mostly tied to non-closed form expression for the posterior covariance matrix, and that the covariance $\hat\Gamma|_{\Theta_c} = \hat\Gamma|_{\Theta_c}(\vecsigma^\ast)$ depends on $\vecsigma$. Active learning type strategies are common in addressing these issues, see~\cite{CABN19, W70, jagalur2021batch, wu2023fast,huan2016sequential,santosa2022bayesian} and we discuss them in \Cref{sec:activeLearning}. 
}

\subsection{Linearization}\label{sec:linearization}
Linearization plays a central role in PDE-based inverse problems. As seen above, derivatives such as $\tfrac{d}{d\vecsigma}\mathrm{pde}$ arise repeatedly in both optimization and Bayesian formulations. Theoretically, linearization is also tied to the notation $B_R(\vecsigma_\ast)$ in the well-posedness estimates~\eqref{eqn:well_posed}, \eqref{eqn:well_posed_semifinite}, and \eqref{eqn:well_posed_finite}. In many cases, well-posedness can  be guaranteed only within a neighborhood of the ground truth $\vecsigma_\ast$. This locality is especially relevant in the finite-dimensional setting: Close to two different underlying truth $\vecsigma_{*,1}$ and $\vecsigma_{*,2}$, the same measurement $\mathrm{d}(\theta_1,\theta_2)$ may exhibit very different sensitivities to parameter variations. Consequently, the type and amount of data needed for stable reconstruction depend on the specific ground truth $\vecsigma_\ast$. 

To address this issue, one must examine how the data $\mathrm{d}(\theta_1,\theta_2)$ varies with perturbations of the parameter $\vecsigma$ near a given $\vecsigma_\ast$. 
Linearization and functional gradient computations provide precisely this information~\cite{stefanov2009linearizing}. 
For PDE-based inverse problems, these tools can be formulated systematically using the calculus of variations. 
Following the notation of~\eqref{eqn:PDE}, the linearized problem typically reduces to a Fredholm integral equation of the first kind:
\begin{equation}\label{eqn:linearization_general}
    \left\langle\frac{d}{d\vecsigma}\mathrm{pde}(\vecsigma,\theta),\delta\vecsigma\right\rangle=\left\langle\frac{d}{d\vecsigma}M_{\theta_2}[u(\vecsigma;\theta_1)],\delta\vecsigma\right\rangle = \delta \mathrm{d}{ {(\theta)}}\,,
\end{equation}
where
\[
\delta\vecsigma=\vecsigma-\vecsigma_\ast \quad\text{and}\quad \delta \mathrm{d} {(\theta)} = \mathrm{d}(\theta)-M_{\theta_2}[u_{\theta_1}(\vecsigma_\ast)]=\mathrm{d}(\theta)-\mathrm{pde}(\vecsigma_*,\theta)\,.
\]
The inverse problem in this linearized form is to use pairs of $\big(\tfrac{d}{d\vecsigma}M_{ {\theta_2}}[u(\vecsigma;{ {\theta_1}})]\,, \delta \mathrm{d} {({\theta})}\big)$, collected over many $\theta$, to reconstruct $\delta\vecsigma$.

While the computation of $\delta \mathrm{d}$ is straightforward, the central challenge lies in evaluating the functional gradient $\tfrac{d}{d\vecsigma}\mathrm{pde}(\vecsigma,\theta) = \tfrac{d}{d\vecsigma}M_{\theta_2}[u(\vecsigma;\theta_1)]$,  {where $u(\vecsigma;\theta_1)$ satisfies} the PDE constraint $\mathcal{L}_{\vecsigma}[u]=S_{\theta_1}$. Different PDEs yield different functional gradients, but a common structure emerges:
\[
\frac{d}{d\vecsigma}M_{\theta_{ {2}}}[u(\vecsigma;\theta_{ {1}})]=f_{\theta_1}\cdot g_{\theta_2}\,
\]
where $f_{\theta_1}$ is associated with the forward PDE solution $u {=u(\sigma,\theta_1)}$ (driven by source $\theta_1$), and $g_{\theta_2}$ is associated with the adjoint PDE solution that is driven by detector $\theta_2$. 
 {The definitions of $f$ and $g$ depend on the specific PDE setting.
With few exceptions, the parameterization $\theta_1$ of the incoming data appears only in the forward equation and thus  in $f$, while information for the detector $\theta_2$ enters only through the adjoint equation, thus appearing only in $g$.  
This separation can be seen \Cref{ex:linearize EIT}.}

 {By} substituting this form into~\eqref{eqn:linearization_general}, renaming $\vecb := \delta \mathrm d$, and denoting $\vecsigma$ as $\delta \vecsigma$, we obtain the following compact expression:
\begin{equation}\label{eq:summary_linearized}
    \left\langle f_{\theta_1}g_{\theta_2},\,  {\delta} \vecsigma \right\rangle = \vecb(\theta_1,\theta_2)\,.
\end{equation}
 {An explicit example is provided below in \Cref{ex:linearize EIT}.} This tensorized structure in the Fredholm integral is a distinctive feature of PDE-based inverse problems, and it will serve as the foundation for the derivations that follow. While linearization itself is often treated as a black-box solver and may be of secondary importance for the purposes of this review, the tensorized structure plays a central role in data selection. 
In particular, effective designs must be compatible with this structure, a requirement that distinguishes PDE-based inverse problems from other experimental design settings.

We now describe the main techniques for computing functional gradients. 
(Readers already familiar with adjoint-based methods may skip this discussion.) 
For simplicity, we fix $\theta=(\theta_1,\theta_2)$ throughout this subsection and omit the subscript.

\medskip
\noindent\textbf{Linearized PDE operator via calculus of variations.}
The general recipe is straightforward. 
Given the PDE constraint $\mathcal{L}_{\vecsigma}[u]=S_{ {\theta_1}}$, we regard $u=u(\vecsigma) {= u(\vecsigma, \theta_1)}$ as a function of the parameter $\vecsigma$. 
Differentiation of the observation functional then yields
\begin{equation}\label{eqn:derivative_sigma}
    \frac{d}{d\vecsigma} M_{ {\theta_2}}[u(\vecsigma)]
    = \left\langle \frac{\delta M_{ {\theta_2}}}{\delta u} \,, \frac{du}{d\vecsigma} \right\rangle  +  {\left.\frac{\delta M_{\theta_2}}{\delta \vecsigma}\right|_{u(\vecsigma)}}.
\end{equation}
Typically, the variational derivative {s} $\frac{\delta M_{ {\theta_2}}}{\delta u}$  {and $\left.\frac{\partial M_{\theta_2}}{\partial \vecsigma}\right|_{u(\vecsigma)}$ are}  explicit; the main difficulty lies in evaluating $\frac{du}{d\vecsigma}$.

To compute this term, we differentiate the PDE constraint. 
Since $\mathcal{L}_{\vecsigma}[u(\vecsigma)]=S_{ {\theta_1}}$ holds identically, we have
\[
\left.\frac{d\mathcal{L}_{\vecsigma}}{d\vecsigma}\right|_{\vecsigma, u(\vecsigma)}
=
\left.\frac{\partial \mathcal{L}_{\vecsigma}}{\partial\vecsigma}\right|_{\vecsigma, u(\vecsigma)}
+ \left.D_u\mathcal{L}_{\vecsigma}\right|_{\vecsigma, u(\vecsigma)} \cdot \frac{du}{d\vecsigma}
=  {\frac{\delta S_{\theta_1}}{\delta \vecsigma}}\,,
\]
or equivalently,
\[
\left.D_u\mathcal{L}_{\vecsigma}\right|_{\vecsigma, u(\vecsigma)} \cdot \frac{du}{d\vecsigma}
= -\left.\frac{\partial \mathcal{L}_{\vecsigma}}{\partial\vecsigma}\right|_{\vecsigma, u(\vecsigma)} +  {\frac{\delta S_{\theta_1}}{\delta \vecsigma}}.
\]
Here  {all three terms } $D_u\mathcal{L}_{\vecsigma}$,  $\tfrac{\partial\mathcal{L}_{\vecsigma}}{\partial\vecsigma}$  { and $\frac{\delta S_{\theta_1}}{\delta \vecsigma}$} are explicit from the PDE. 
Direct inversion of $D_u\mathcal{L}_{\vecsigma}$, however, is often prohibitively expensive.

A standard alternative is the \textbf{adjoint method}. Returning to~\eqref{eqn:derivative_sigma}, we introduce an adjoint variable $\lambda$ defined by
\begin{equation}\label{eqn:adjoint}
    \left(D_u\mathcal{L}_{\vecsigma}\right)^\ast \lambda
    = \frac{\delta M_{ {\theta_2}}}{\delta u}\,.
\end{equation}
Multiplying the differentiated PDE constraint by $\lambda$ and integrating by parts yields:
\begin{equation}\label{eqn:functional_derivative_final}
\begin{aligned}
    \frac{d}{d\vecsigma} M_{ {\theta_2}}[u(\vecsigma)] -    {\left.\frac{\delta M_{\theta_2}}{\delta \vecsigma}\right|_{u(\vecsigma)}}
    &= \left\langle \frac{\delta M_{ {\theta_2}}}{\delta u}, \frac{du}{d\vecsigma} \right\rangle
    = \left\langle \left(D_u\mathcal{L}_{\vecsigma}\right)^\ast \lambda, \frac{du}{d\vecsigma} \right\rangle \\
    &= \left\langle \lambda, D_u\mathcal{L}_{\vecsigma} \cdot \frac{du}{d\vecsigma} \right\rangle
    = -\left\langle \lambda, \left.\frac{\partial \mathcal{L}_{\vecsigma}}{\partial\vecsigma}\right|_{\vecsigma, u(\vecsigma)}   {-\frac{\delta S_{\theta_1}}{\delta \vecsigma}}\right\rangle.
\end{aligned}
\end{equation}
This adjoint formulation is highly efficient: only one forward PDE solve and one adjoint PDE solve \eqref{eqn:adjoint} are required, independent of the dimension of $\vecsigma$.

Since $M_{ {\theta_2}}[u(\vecsigma)]$ maps a function $\vecsigma$ to a scalar, its derivative with respect to $\vecsigma$ is itself a function. 
Consequently, the terms $\frac{du}{d\vecsigma}$ and $\frac{\partial \mathcal{L}_{\vecsigma}}{\partial \vecsigma}$ in~\eqref{eqn:functional_derivative_final} should be interpreted as operators, with the inner product taken in the operator–function sense. 
This viewpoint naturally motivates a projected formulation, which is often more convenient in practice.
We have
\begin{equation}\label{eqn:functional_proj}
\Big\langle \tfrac{d}{d\vecsigma} M_{ {\theta_2}}[u(\vecsigma)], \delta\vecsigma \Big\rangle
= -\left\langle \lambda, \tfrac{\partial\mathcal{L}_{\vecsigma}}{\partial\vecsigma}[\delta\vecsigma]  {- \tfrac{\delta S_{\theta_1}}{\delta \vecsigma}[\delta \vecsigma] }\right\rangle   {+\left\langle\left.\tfrac{\delta M_{\theta_2}}{\delta \vecsigma}\right|_{u(\vecsigma)},\delta \vecsigma\right \rangle},
\end{equation}
where $\tfrac{\partial\mathcal{L}_{\vecsigma}}{\partial\vecsigma}[\delta\vecsigma]$ denotes the action of $\tfrac{\partial\mathcal{L}_{\vecsigma}}{\partial\vecsigma}$ on $\delta\vecsigma$. 
This projected formulation avoids explicit construction of the full functional derivative and directly yields directional gradients, making it highly practical.

\begin{example}[Linearized EIT \ref{itm: E2}]\label{ex:linearize EIT}
 {To cast~\eqref{eqn:elliptic}--\eqref{eqn:elliptic:DirichletBoundary} into the form $\mathcal L_{\vecsigma} [u] = S$,  we homogenize the boundary. Choosing  $\tilde \phi\in H^1(\Omega)$ as an extension of $\phi$, i.e. with $\tilde \phi\mid_{\partial \Omega} = \phi$,  allows us to rewrite
~\eqref{eqn:elliptic}--\eqref{eqn:elliptic:DirichletBoundary}  equivalently  in terms of $w = u - \tilde \phi$ as 
\begin{equation}
    \nabla \cdot \big(\vecsigma(x) \nabla w\big) = - \nabla \cdot \big(\vecsigma(x) \nabla \tilde \phi\big) 
\end{equation}
with homogeneous boundary condition $w\mid_{\partial \Omega} = 0$. This equation has the desired form, with PDE operator defined through $\mathcal L_{\vecsigma}[w]:= \nabla \cdot \big(\vecsigma(x) \nabla w\big)$, and the data $\phi$ determines the right hand side $S[\phi] = - \nabla \cdot \big(\vecsigma(x) \nabla \tilde \phi\big)$. The measurement becomes $ M_{\theta_2}[w(\vecsigma)]= M_{\theta_2}[u(\vecsigma)] = \vecsigma \partial_{n}(w+\tilde \phi)(\theta_2)$.

The ingredients of the functional derivative are as follows.

\begin{enumerate}
    \item \textbf{Operator derivative} $\frac{\partial\mathcal{L}_{\vecsigma}}{\partial\vecsigma}$.  
    By definition, we have
    \[
    \left.\frac{\partial\mathcal{L}_{\vecsigma}}{\partial\vecsigma}\right|_{\vecsigma,w}[\delta\vecsigma]
    = \left.\frac{d}{d\epsilon} \mathcal{L}_{\vecsigma+\epsilon\delta\vecsigma}[w]\right|_{\epsilon=0},
    \]
    where the expression is evaluated at $(\vecsigma, w)$ and applied in the direction $\delta\vecsigma$.  
    Observing that
    \[
    \mathcal{L}_{\vecsigma+\epsilon\delta\vecsigma}[w]
    = \nabla\cdot\big((\vecsigma+\epsilon\delta\vecsigma)\nabla w\big)
    = \nabla\cdot(\vecsigma\nabla w) + \epsilon\,\nabla\cdot(\delta\vecsigma\nabla w),
    \]
    we obtain
    \[
    \left.\frac{\partial\mathcal{L}_{\vecsigma}}{\partial\vecsigma}\right|_{\vecsigma,w}[\delta\vecsigma]
    = \nabla\cdot(\delta\vecsigma\nabla w).
    \]
    
    \item \textbf{Fr\'echet derivative in $w$}, $D_w\mathcal{L}_{\vecsigma}$.  
    We know that
    \[
    D_w\mathcal{L}_{\vecsigma}[\delta w]
    = \nabla\cdot(\vecsigma\nabla \delta w).
    \]
    This operator is self-adjoint, so we have
    $\left(D_w\mathcal{L}_{\vecsigma}\right)^\ast \lambda
    = \nabla\cdot(\vecsigma\nabla \lambda)$. 
    \item  {\textbf{Measurement derivative.} The Fr\'echet derivative of the Neumann measurement in $\vecsigma$ reads $\left.\tfrac{\delta M_{\theta_2}}{\delta \vecsigma}\right|_{w(\vecsigma)} = \partial_n (w(\theta_2)+ \tilde \phi)$.  }
    \item  {\textbf{Input derivative.} The Fr\'echet derivative  of the input data w.r.t. $\vecsigma$ can be computed equivalently to the operator derivative as $\left.\tfrac{\delta S_{\theta_1}}{\delta \vecsigma}\right.[\delta \vecsigma] =- \nabla \cdot(\delta \vecsigma \nabla \tilde \phi)$.  }
\end{enumerate}
By combining these results with~\eqref{eqn:functional_proj}, we obtain
\begin{align*}
 \left\langle\frac{d}{d\vecsigma}M_{ {\theta_2}}[w(\vecsigma)]\,,\delta\vecsigma\right\rangle
&= -\left\langle \lambda,\, \nabla\cdot(\delta\vecsigma\nabla (w+\tilde \phi) \right\rangle  {+ (\delta \vecsigma)\partial_n(w+\tilde \phi)(\theta_2)} \\
&= \left\langle \nabla\lambda \cdot \nabla (w+\tilde \phi),\, \delta\vecsigma \right\rangle   {+ \langle \delta \vecsigma\partial_n(w+\tilde \phi), \delta_{\theta_2} -\lambda\rangle_{\partial \Omega}}.
\end{align*}
where we used integration-by-parts, and formally denoted by $\langle \cdot ,\delta_{\theta_2} -\lambda\rangle_{\partial \Omega}$  the respective functional on boundary data.  {By picking a suitable Dirichlet boundary condition for $\lambda$, namely $\lambda = \delta_{\theta_2}$, the second term vanishes.} We transform back to $u$ and  reinstate the dependence on $\theta=(\theta_1,\theta_2)$  to obtain
\begin{equation}\label{eq:linear eit cont}
\tfrac{d}{d\vecsigma}\mathrm{pde}(\vecsigma,\theta)
= \nabla \lambda_{\theta_2} \cdot \nabla u_{\theta_1} ,
\end{equation}
where $u_{\theta_1}$ and $\lambda_{\theta_2}$ are solutions of the forward and adjoint problems, respectively. The linearized inverse problem becomes
\begin{equation}\label{eq:EIT linearized}
\Big\langle \nabla\lambda_{\theta_2} \cdot \nabla u(\vecsigma^\ast;\theta_1), \delta\vecsigma \Big\rangle
= \mathrm{d}(\theta) - M_{ {\theta_2}}[u(\vecsigma^\ast;\theta_1)]\,.
\end{equation}
}
\end{example}

This example highlights the decoupling of $u$ and $\lambda$ associated with source $\theta_1$ and detector $\theta_2$, respectively, in consistency with the tensorized structure discussed in~\eqref{eq:summary_linearized}.

\section{Randomized linear algebra techniques}\label{sec:rnla}
Randomized numerical linear algebra (RNLA) incorporates randomization into the design of algorithms for numerical linear algebra algorithms.
In contrast to classical deterministic methods, randomized algorithms deliberately trade accuracy for efficiency. 
They are particularly effective for large-scale problems or for settings that require repeated linear operations, where computational cost rather than a need for high precision is the primary bottleneck. 
Early applications of RNLA arose in such data-science domains as imaging, signal processing, and internet-scale problems, where speed is of paramount importance~\cite{martinsson2020randomized, mahoney2011randomized, W14}. 
More recently, applications in the physical sciences have emerged, where accuracy demands must be balanced against computational feasibility, making RNLA both powerful and practical.

PDE-based experimental design aligns naturally with this perspective. 
Inverse problems governed by PDEs ultimately prioritize reconstruction accuracy. 
Yet before reconstruction can take place, one must often examine a large pool of potential data to identify informative subsets. 
This selection step is typically repeated many times and, crucially, it does not require exact precision. 
RNLA algorithms are well-suited to this step, as they accelerate the repeated linear algebra computations underlying data selection while maintaining enough accuracy to guide the process.

A classical entry point for RNLA in experimental design is through leverage scores. 
Leverage scores are statistical quantities that measure the importance of each constraint in a linear regression problem, and they have long been used to guide data selection. 
However, they are expensive to compute exactly, which is why randomization becomes valuable. 
To illustrate, consider the linear regression problem:
\[
\min_{\vecsigma}\|\mA\vecsigma-\vecb\|_2\,.
\]
Here $\mA\in\R^{C\times n}$ is the design matrix, with each row encoding one constraint, and $\vecb$ collects the data.  {When $\mA$ has full column rank, the} optimal solution and its projection are explicitly given by
\[
\vecsigma^*=(\mA^\top\mA)^{-1}\mA^\top\vecb\,,\quad\text{and}\quad\vecb^*=\mA(\mA^\top \mA)^{-1}\mA^\top \vecb\,.
\]
The influence of the data $\vecb$ on the fitted output $\vecb^*$ can be quantified through the Jacobian:
\[
\nabla_{\vecb}\vecb^*=\mA(\mA^\top\mA)^{-1}\mA^\top\doteq \mathbf{H}\,,
\]
known as the hat matrix (as it is often, though not here, denoted by $\hat{\mathbf{H}}$). 
Its diagonal entries, $l_i=\mathbf{H}_{ii}$ define the leverage scores, which measure the relative importance of the $i$-th observation. 
High-leverage points correspond to experiments with strong influence on the regression fit. 
Leverage scores are widely used, for example, in outlier detection in economics and medical imaging ~\cite{leverage_economics, mejia2017pca}, kernel methods \cite{avron17a, lee_kernel} and graph learning \cite{spielman2008graph}. 

Computing $\mathbf{H}$ presents challenges in computation and data access. 
Forming this matrix directly requires $O(Cn^2)$ flops, which is costly for large problems. 
Even more problematic is the sample complexity: Evaluating leverage scores requires explicit knowledge of the entire matrix $\mA$, meaning all experiments must be conducted and formulated before any selection is possible. 
Randomized methods offer a way to address both challenges: They reduce the computational cost to near-linearity in $C$ and $n$ and allow approximate leverage scores to be estimated from only partial or sampled information about $\mA$~\cite{liberty2007randomized,DMMW12,martinsson2020randomized}.
 
While these techniques were not originally developed for PDE-based inverse problems, they represent one of the earliest demonstrations of how RNLA can be useful for experimental design. 
Over the past decade, this perspective has broadened: Methods such as matrix sketching, randomized matrix–matrix multiplication, column selection, randomized singular value decomposition (SVD), compressed sensing, and matrix completion have all been adapted to PDE settings. 
Collectively, they provide ways to accelerate and guide data selection. 
In what follows, we first review these methods in their general linear algebraic form, and then describe how they can be specialized for PDE-based experimental design.

\subsection{Matrix sketching}
Matrix sketching is an important problem class in RNLA. 
Conceptually, it can be understood in analogy to regression: in classical least-squares regression, one may ``sketch" by reducing the row dimension of the system in a principled way, hoping to obtain a solution that is nearly as good as the full system.

Formally, consider the overdetermined linear system
\[
\mA\, \vecsigma \approx \vecb\,, \quad \mA\in \R^{C\times n}\,,C\gg n\,,
\]
whose least-squares solution is
\begin{equation}\label{eq:linear ls}
\vecsigma^* = \arg\min_{\vecsigma \in \R^n} \| \mA \vecsigma - \vecb \|_2^2 \,.
\end{equation}
Since the system is highly overdetermined, a reasonable approximation can often be obtained by using only a reduced set of rows from $\mA$.
To achieve this, one introduces a sketching matrix $\mS \in \R^{c \times C}$ with $c \ll C$ and solves the reduced problem:
\begin{equation}\label{eq:sketched linear}
\vecsigma^*_{\mS} \ = \ \arg\min_{\vecsigma \in \R^n} \| \mS\mA \vecsigma - \mS\vecb \|_2^2\,,
\end{equation}
with the goal that $\vecsigma^*_{\mS} \approx \vecsigma^*$. A typical guarantee is:
\begin{equation}\label{eqn:sketch_goal}
\| \mA {\vecsigma}^*_\mS -\vecb \|_2^2 \leq (1+\epsilon)\| \mA \vecsigma^* - \vecb \|_2^2, \quad \text{with~probability~}1-\delta.
\end{equation}
The mechanism is illustrated in Figure~\ref{fig: sketching}.
\begin{center}
\resizebox{0.85\columnwidth}{!}
{\begin{tikzpicture}

\draw[fill=blue!10, draw=blue, thick] (-0.7,0.2) rectangle (1.8,0.8);
\node[rectangle,minimum width=0cm, minimum height = 0cm, align=center,fill=white,fill opacity=0,text opacity=1,font=\bfseries,label=right:{$\bf S$}] at (1,0.4) (block) {} ; 
\node[rectangle,minimum width=0cm, minimum height = 0cm, align=center,fill=white,fill opacity=0,text opacity=1,font=\bfseries,label=right:{$\boldsymbol{\cdot}$}] at (1.7,0.5)  {} ; 
\draw [decoration={brace,amplitude=4pt,mirror,},decorate,line width=1pt] ($(block)+(-4.5em,2ex)$) -- ($(block)+(-4.5em,-2ex)$);
\draw [decoration={brace,amplitude=4pt},decorate,line width=1pt] ($(block)+(3.5em,-9ex)$) -- ($(block)+(3.5em,10ex)$);
\node[rectangle,minimum width=0cm, minimum height = 0cm, align=center,fill=white,fill opacity=0,text opacity=1,font=\bfseries,label=right:{$c$}] at (-1.7,0.4)  {} ; 
\node[rectangle,minimum width=0cm, minimum height = 0cm, align=center,fill=white,fill opacity=0,text opacity=1,font=\bfseries,label=right:{$C$}] at (1.7,1)  {} ;

\draw[fill=red!10, draw=red, ultra thick] (2.5,-1) rectangle (3.2,2);
\node[rectangle,minimum width=0cm, minimum height = 0cm, align=center,fill=white,fill opacity=0,text opacity=1,font=\bfseries,label=right:{$\bf A$}] at (2.4,0.4) (block2) {} ; 
 
\draw [decoration={brace,amplitude=4pt,mirror,},decorate,line width=1pt] ($(block2)+(2em,10ex)$) -- ($(block2)+(0.2em,10ex)$);
\node[rectangle,minimum width=0cm, minimum height = 0cm, align=center,fill=white,fill opacity=0,text opacity=1,font=\bfseries,label=right:{$n$}] at (2.5,2.45)  {} ; 

\draw[draw=black, line width=3pt, rounded corners=1pt] (3.6,0.2) -- (3.6,0.8);
\node[rectangle,minimum width=0cm, minimum height = 0cm, align=center,fill=white,fill opacity=0,text opacity=1,font=\bfseries,label=right:{\large$\hat{\pmb{\sigma}}$}] at (3.5,0.5) {} ;

\draw[thick] (4.2,0.6) -- (4.5,0.6);
\draw[thick] (4.2,0.4) -- (4.5,0.4);

\draw[fill=blue!10, draw=blue, thick] (5,0.2) rectangle (7.5,0.8);
\node[rectangle,minimum width=0cm, minimum height = 0cm, align=center,fill=white,fill opacity=0,text opacity=1,font=\bfseries,label=right:{$\bf S$}] at (6.5,0.4)  {} ; 

\node[rectangle,minimum width=0cm, minimum height = 0cm, align=center,fill=white,fill opacity=0,text opacity=1,font=\bfseries,label=right:{$\boldsymbol{\cdot}$}] at (7.4,0.5)  {} ; 

\draw[draw=red, ultra thick] (8,-1) -- (8,2);
\node[rectangle,minimum width=0cm, minimum height = 0cm, align=center,fill=white,fill opacity=0,text opacity=1,font=\bfseries,label=right:{$\bf b$}] at (7.8,0.4)  {} ; 


\node[rectangle,minimum width=0cm, minimum height = 0cm, align=center,fill=white,fill opacity=0,text opacity=1,font=\bfseries,label=right:{$\Longrightarrow$}] at (8.8,0.4)  {} ; 

\draw[fill=orange!30, draw=orange, thick] (11.5,0.2) rectangle (12.3,0.8);
\node[rectangle,minimum width=0cm, minimum height = 0cm, align=center,fill=white,fill opacity=0,text opacity=1,font=\bfseries,label=right:{$\bf S A$}] at (11.4,0.4) (block3) {} ; 

\draw [decoration={brace,amplitude=4pt,mirror,},decorate,line width=1pt] ($(block3)+(0em,2ex)$) -- ($(block3)+(0em,-2ex)$);
\draw [decoration={brace,amplitude=4pt,mirror,},decorate,line width=1pt] ($(block3)+(2.2em,3ex)$) -- ($(block3)+(0em,3ex)$);

\node[rectangle,minimum width=0cm, minimum height = 0cm, align=center,fill=white,fill opacity=0,text opacity=1,font=\bfseries,label=right:{$c$}] at (10.5,0.5)  {} ; 

\node[rectangle,minimum width=0cm, minimum height = 0cm, align=center,fill=white,fill opacity=0,text opacity=1,font=\bfseries,label=right:{$n$}] at (11.5,1.3)  {} ;

\draw[draw=black, line width=3pt, rounded corners=1pt] (12.8,0.2) -- (12.8,0.8);

\node[rectangle,minimum width=0cm, minimum height = 0cm, align=center,fill=white,fill opacity=0,text opacity=1,font=\bfseries,label=right:{\large$\hat{\pmb{\sigma}}$}] at (12.7,0.5) {} ;
\draw[thick] (13.4,0.6) -- (13.7,0.6);
\draw[thick] (13.4,0.4) -- (13.7,0.4);

\draw[fill=orange, draw=orange, ultra thick] (14,0.2) rectangle (14,0.8);
\node[rectangle,minimum width=0cm, minimum height = 0cm, align=center,fill=white,fill opacity=0,text opacity=1,font=\bfseries,label=right:{$\bf S b$}] at (13.8,0.4)  {} ; 

\end{tikzpicture}}

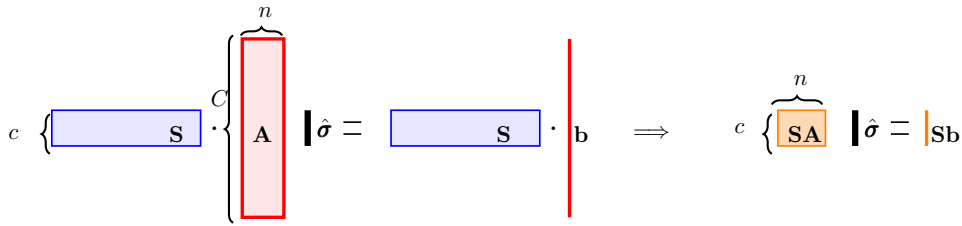
\captionof{figure}{Sketching mechanism.}
\label{fig: sketching}
\end{center}

The advantages of sketching are twofold. 
From a \textbf{computational perspective}, solving~\eqref{eq:linear ls} requires $O(C n^2)$ flops, while~\eqref{eq:sketched linear} requires only $O(\min(cn^2, c^2 n))$ with $c \ll C$, yielding significant savings. From an \textbf{experimental design} perspective, the savings are even more pronounced in sample complexity: the full problem requires $C$ measurements (the rows of $\mA$ and entries of $\vecb$), whereas the sketched problem reduces this to $c$ measurements, substantially cutting the effort in data acquisition.

\begin{figure}[!htb]
\centering
\subfloat[subGaussian]{\resizebox{0.2\textwidth}{!}{\begin{tikzpicture}
  \matrix[matrix of nodes,
          nodes in empty cells,
          nodes={minimum size=10mm, outer sep=1pt, draw = none, fill=orange},
          column sep=2mm, row sep=2mm,
          draw] (m) {
    \node {}; & \node {}; & \node(target) {}; & \node {}; & \node {};  \\
    \node {}; & \node {}; & \node {}; & \node {}; & \node {}; \\
  };

  \node[above right=0.5cm and 0.5cm of target] (label) {$\Large{\sim \mathcal{N}(0,\frac{1}{c})}$};
  \draw[->, thick] (label.south west) -- (target.north east);
\end{tikzpicture}}}~~~~
\subfloat[FJLT]{\resizebox{0.27\textwidth}{!}{\begin{tikzpicture}
  \matrix[matrix of nodes,
          nodes in empty cells,
          nodes={minimum size=10mm, outer sep=1pt, draw = none, fill=orange},
          column sep=2mm, row sep=2mm,
          draw] (m) {
            \node(t1) {} ; & \node(t2) {}; & \node(t3) {}; & \node(t4) {}; &\node(t5) {}; \\
    \node[fill=gray!50] () {}; & \node[fill=gray!50] {}; & \node[fill=gray!50] {}; & \node[fill=gray!50] {}; & \node[fill=gray!50] {}; \\
    \node[fill=gray!50] {}; & \node[fill=gray!50] {}; & \node[fill=gray!50] {}; & \node[fill=gray!50] {}; & \node[fill=gray!50] {}; \\
      \node(t21) {}; & \node {}; & \node {}; & \node {}; &\node {};  \\
        \node[fill=gray!50] {}; & \node[fill=gray!50] {}; & \node[fill=gray!50] {}; & \node[fill=gray!50] {}; &\node[fill=gray!50] {};  \\      
  };

     \node[left=0.5cm of t21] (label) {\Large${\pm\frac{1}{\sqrt{n}} \omega^{jk}}$};
  \draw[->, thick] (t21.west) -- (label.east);
\end{tikzpicture}}}~~~~
\subfloat[sparse]{\resizebox{0.2\textwidth}{!}{\begin{tikzpicture}
  \matrix[matrix of nodes,
          nodes in empty cells,
          nodes={minimum size=10mm, outer sep=1pt, draw = none, fill=gray!50},
          column sep=2mm, row sep=2mm,
          draw] (m) {
    \node {}; & \node {}; & \node[fill=orange] (target) {}; & \node {}; & \node {};  \\
    \node {}; & \node {}; & \node {}; & \node {}; & \node[fill=orange] {}; \\
  };

  \node[above right=0.5cm and 0.5cm of target] (label) {$\Large{\in \{+1 ,-1 \}}$};
  \draw[->, thick] (label.south west) -- (target.north east);

\end{tikzpicture}}}
\caption{Major classes of sketching matrix choices. Orange squares indicate nontrivial entries. Gray squares represent unsampled entries in (b), and zero-valued entries in (c). In (b) FJLT, the $(j,k)$-th entry of $\pmb{\mathcal{F}}$ is given by $\pm \frac{1}{\sqrt{n}} \omega^{jk}$, where $\omega = \text{exp}(-\frac{2 \pi \mathrm{i}}{n})$ is the $n$-th root of unity, and $jk$ is the exponent. The sign (positive or negative) is determined by a random sign-flip applied to its column.}\label{fig:design_principle_sketching}
\end{figure}
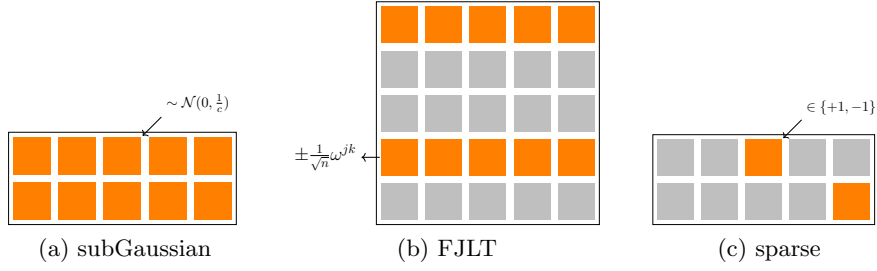

The main challenge is to construct a sketching matrix $\mS$ that ensures the guarantee~\eqref{eqn:sketch_goal} while maintaining the property $c \ll C$. 
Without entering technical details, we highlight several widely used strategies (see~\cite{W14, M16}):
\begin{itemize}
    \item SubGaussian matrices \cite{RV08}.     
    The sketching matrix is generated by i.i.d. zero-mean, isotropic subGaussian entries. Common subGaussian types include Gaussian, Bernoulli, uniform distributions. 
    \item Fast Johnson-Lindenstrauss transform (FJLT) \cite{AC09,tropp2011improved}.    
    The sketching matrix is constructed as $\pmb{\Pi}\pmb{\mathcal{F}},$ where $\pmb{\mathcal{F}} \in \R^{C \times C}$ represents a discrete Fourier or Hadamard transform with randomly sign-flipped columns  and $\pmb{\Pi} \in \R^{c \times C}$ denotes a random row sampling operation.
    
    \item Sparse sketching matrices \cite{A03, CHARIKAR20043}. Particularly useful in streaming settings, these matrices have entries in \(\{0, +1, -1\}\), chosen according to a suitable distribution. Examples include CountSketch~\cite{CHARIKAR20043} and sparse random projection~\cite{A03}.  
\end{itemize}
The three design principles are visualized in Figure~\ref{fig:design_principle_sketching}.

While the precise dependence of $c$ on $\epsilon$ and $\delta$ varies across these constructions, their theoretical justification is unified by the celebrated Johnson–Lindenstrauss lemma~\cite{JL84}, which guarantees the existence of low-dimensional embeddings that preserve geometry with high probability.

\subsubsection{Matrix sketching for PDE based experimental design} 
Sketching aims to provide an almost equally good reconstruction to a regression problem using only a subsampled system. This perspective aligns naturally with experimental design. However, PDE-based inverse problems exhibit a special structure that prevents a direct use of classical sketching. Most notably, as shown in~\eqref{eq:summary_linearized}, the linearized PDE-based inverse problem typically takes a tensorized form.

 We rewrite the linearized system in algebraic form:
\begin{equation}\label{eq:eit system}
{\bf A} \cdot\vecsigma = {\bf b} \,,
\end{equation}
where $\bf b$ encodes the right-hand side of~\eqref{eqn:linearization_general} (with entries indexed by $\theta$), and $\bf A$ is essentially the Jacobian matrix.  {In the interest of readability, here and in the following, we switch the notation of the parameter in the linearized problem from $\delta\vecsigma$ to $\vecsigma$, assuming that the considered version is clear from the context.}

Unlike in classical regression, $\bf A$ is infinite in size: its rows are indexed by $\theta=(\theta_1,\theta_2)\in\Theta$, the entire design space, while its columns are indexed by the discretization of $\vecsigma$. 
Restricting to $\vecsigma\in\Sigma_n$ with dimension $n$, each row of $\mA_{\theta,:}$ lies in $\R^n$. 

The key distinction from standard regression is the tensorized structure of $\mA$. In numerical setting of sketching, we focus on the finite but possibly large data space as $\Theta_C$. By~ {\eqref{eq:summary_linearized}}, $\mA$ can be expressed columnwise product.
Let
\[
\mathbf{F}\in\R^{C_1\times n}\,, \quad \mathbf{G}\in\R^{C_2\times n}\quad\text{with}\quad\mathbf{F}_{\theta_1,:} = f_{\theta_1}\,,\quad\mathbf{G}_{\theta_2,:} = g_{\theta_2}\,,
\]
where $\theta_i\in\Theta_i$ and $|\Theta_i|=C_i$ and  { $\Theta = \Theta_1\times \Theta_2$ with} $C=C_1C_2$. Then 
\[
\mA = \mathbf{F}\ast\mathbf{G}\in\R^{C_1C_2\times n}\quad\text{with each row}\quad \mA_{\theta,:} = \mathbf{F}_{\theta_1,:}\  {\odot}\ \mathbf{G}_{\theta_2,:} \in \mathbb{R}^n  \qquad  { \text{for } \theta = (\theta_1,\theta_2)\in \Theta}\,,
\]
 {where $\odot$ denotes the entry-wise product, i.e. each entry in the row vector ${\bf A}_{\theta,i} = {\bf F}_{\theta_1,i} {\bf G}_{\theta_2,i},$ for $i = 1, \dots, n,$ and $\ast$ is the Khatri-Rao tensor product. The measurement index $\theta$ encapsulates two indices: source $\theta_1$ and detection location $\theta_2.$}
This structure requires a tensor-aware sketching procedure, with the sketching matrix also having a row-wise tensorized form, so that $\mathbf{F}$ and $\mathbf{G}$ are sketched independently. Physically, this corresponds to the fact that the parameter specifying the PDE input ($\theta_1$) and the parameter specifying the detector ($\theta_2$) are independently configurable in a fully factorized design.
Two representative constructions along these lines are the following.
\begin{itemize}
    \item SubGaussian tensorizing. Each row of the sketching matrix $\mS$ is written as $\mS_{i,:} = \mS^{\mathbf{F}}_{i,:}\otimes \mS^{\mathbf{G}}_{i,:}$, where $\mS^{\mathbf{F}}$ and $\mS^{\mathbf{G}}$ are filled with i.i.d. zero-mean isotropic subGaussian entries.
    \item Kronecker FJLT. The transform $\pmb{\Pi}\pmb{\mathcal{F}}$ is modified to the tensorized form $\pmb{\Pi}\big(\pmb{\mathcal{F}}^{\mathbf{F}}\otimes\pmb{\mathcal{F}}^{\mathbf{G}}\big)$, so that $\pmb{\mathcal{F}}^{\mathbf{F}}$ and $\pmb{\mathcal{F}}^{\mathbf{G}}$ act separately on $\mathbf{F}$ and $\mathbf{G}$.
\end{itemize}
These mechanisms are illustrated in Figure~\ref{fig: tensor sketching}.
\begin{center}
\resizebox{\columnwidth}{!}
{\begin{tikzpicture}[
 line/.style ={draw, thick, -latex}                ]
 \tikzset{arr/.style={thick,->,>={Stealth[scale=1.5]}}};
 \node[rectangle,minimum width=3cm, minimum height = 0.15cm, align=center,fill=blue,fill opacity=0.4,text opacity=1,font=\bfseries,label=above:{$\bf S$}] at (-6,1) {} ;         
\node[rectangle,minimum width=0.15cm, minimum height = 3cm, align=center,fill=red,fill opacity=0.4,text opacity=1,font=\bfseries,label=right:{$\bf A$}] at (-4,0)  {} ; 

\draw [arr](-2.5,0) --++ (1,0);

\node[rectangle,minimum width=1cm, minimum height = 0.15cm, align=center,fill=blue,fill opacity=0.4,text opacity=1,font=\bfseries,label=above:{${\bf S}^{\mathbf{F}}$}] at (0,1) {} ;  
\node[rectangle,minimum width=0.15cm, minimum height = 1cm, align=center,fill=red,fill opacity=0.4,text opacity=1,font=\bfseries,label=right:{$\mathbf{F}$}] at (1,1) {} ;

\node[rectangle,minimum width=1cm, minimum height = 0.15cm, align=center,fill=blue,fill opacity=0.4,text opacity=1,font=\bfseries,label=above:{${\bf S}^{\mathbf{G}}$}] at (0,-1) {} ;
\node[rectangle,minimum width=0.15cm, minimum height = 1cm, align=center,fill=red,fill opacity=0.4,text opacity=1,font=\bfseries,label=right:{$\mathbf{G}$}] at (1,-1) {} ;

\draw [dashed, thick] (3,1.5) -- (3,-2);

\node[rectangle,minimum width=3cm, minimum height = 3cm, align=center,fill=blue,fill opacity=0.4,text opacity=1,font=\bfseries] at (6,0) (block) {$\pmb{\mathcal{F}}$} ;                
\node[rectangle,minimum width=0.15cm, minimum height = 3cm, align=center,fill=red,fill opacity=0.4,text opacity=1,font=\bfseries,label=right:{$\bf A$}] at (8,0) {} ;

\node[rectangle,minimum width=1cm, minimum height = 1cm, align=center,fill=blue,fill opacity=0.4,text opacity=1,font=\bfseries, label=above:{$\pmb{\mathcal{F}}^{\mathbf{F}} $}] at (12,1) {} ;  
\node[rectangle,minimum width=0.15cm, minimum height = 1cm, align=center,fill=red,fill opacity=0.4,text opacity=1,font=\bfseries,label=right:{$\mathbf{F}$}] at (13,1) {} ;

\node[rectangle,minimum width=1cm, minimum height = 1cm, align=center,fill=blue,fill opacity=0.4,text opacity=1,font=\bfseries,label=above:{$\pmb{\mathcal{F}}^{\mathbf{G}}$}] at (12,-1) {} ;  
\node[rectangle,minimum width=0.15cm, minimum height = 1cm, align=center,fill=red,fill opacity=0.4,text opacity=1,font=\bfseries,label=right:{$\mathbf{G}$}] at (13,-1)  {} ;

\draw [arr](9.2,0) --++ (1,0);

\end{tikzpicture}}
\captionof{figure}{Tensor sketching acceleration: the complexity of computing $\mS\mA$ breaks down to the scale of each factor matrix size.}
\label{fig: tensor sketching}
\end{center}

Tensor-structured sketching has been extensively developed in recent years~\cite{MSW19,CLNW20,JKW20,CJ21,INRZ21}, often relying on the distributive property $(\bf B\otimes \bf D)(C\otimes E) = (BC)\otimes (DE)$. 
In PDE-based inverse problems, the Khatri–Rao product induces the tensorized requirement, with notable works including~\cite{CLNW20, CJ21}. 
Likewise,~\cite{JKW20} developed Kronecker Fourier-based random projections (Kronecker FJLT), extending standard FJLT.

 {In the context of data selection for PDE-based inverse problem, $\mathbf F$ and $\mathbf G$ encode source and detector information respectively, and the sketching has to be conducted independently to ensure a physically meaningful experimental setup. While sparse sketching explicitly selects a small number of representative experimental configurations, the Gaussian tensoring is to perform linear combination of input data and measurements. For example, if the forward  map  $f_{\theta_1}$ is linear with respect to the input data $S_{\theta_1}$, as in \ref{itm: E2}--\ref{itm: E4}, the sketched rows of $\mathbf F$ corresponds to setting up an experiment with the respective weights that linearly combine the input data:
$$(\mS^{\mathbf F}\mathbf F)_{i,:} = \sum_{\theta_1} \mS^{\mathbf F}_{i,\theta_1}\mathbf F_{\theta_1,:} =  \sum_{\theta_1} \mS^{\mathbf F}_{i,\theta_1}f_{\theta_1,:} = f_{\hat \theta_1,:}, $$
where $f_{\hat \theta_1,:}$ corresponds to input  data $S_{\hat \theta_1} = \sum_{\theta_1} \mS^{\mathbf F}_{i,\theta_1}S_{ \theta_1}$.}

Finally, we highlight the issue of sampling complexity: How many  randomly sampled measurements are sufficient to ensure that the sketched solution $\vecsigma^*_\mS$ approximates the full solution? 
A trade-off emerges: The tensor structure accelerates computations but introduces higher-order Gaussian chaos, complicating probabilistic analysis. 
For instance, Theorem 1.1 in~\cite{CJ21} establishes a nearly optimal bound:
\begin{equation} \label{eqn: tensor sample} c \geq \max \left(\mathcal{O} \left(\frac{\text{rank}^2({\bf A}) + \log^2 (1/\delta)}{\epsilon} \right), \mathcal{O} \left(\frac{\text{rank}({\bf A}) + \log (1/\delta)}{\epsilon^2} \right)\right),\end{equation}
which ensures that $\vecsigma^*_\mS$ satisfies the guarantee~\eqref{eqn:sketch_goal} with high probability. 
The tensor sample complexity estimate \eqref{eqn: tensor sample} yields a slightly weaker bound compared to the optimal result regarding standard unstructured sketching \cite{pilanci2015randomized}, given by $c \geq \mathcal{O} (\frac{\text{rank}(\mA) + \log (1/\delta)}{\epsilon^2})$, which corresponds to the second term in \eqref{eqn: tensor sample}.

Due to the inherent source-detector structure in data arising from PDE inverse problems, tensor sketching is specially designed to efficiently process such structured data. Although great progress has been made in studying tensor sketching methods, several key open questions remain.
A central question is to establish the optimal sample complexity bound, improving upon the sub-optimal result in \eqref{eqn: tensor sample}. 
Another crucial question concerns the practical scenario in which measurements
are sparse or missing. 
In such cases, it may be beneficial to design the sketching operation in a correspondingly sparse format, allowing the use of even fewer measurements while still preserving the essential information of the inference problem. 

\subsection{Matrix-Matrix product}\label{sec: sampling}
Matrix multiplication is one of the most fundamental tasks in numerical linear algebra. Given two matrices $\mA \in \R^{m \times C}$ and $\mathbf{B} \in \R^{C \times n}$, we want to compute their product
\[
\mC=\mA\mB\in\R^{m\times n}\,,
\]
for which the formula is $\mC_{ij} = \mA_{i,:}\mB_{:,j} = \sum_k\mA_{ik}\mB_{kj}$ for every entry $(i,j)$. 
Randomization enters through an alternative but equivalent representation of the product:
\begin{equation}\label{eqn:matrix_matrix_product}
\mC=\sum_{k=1}^C\mA_{:,k}\mB_{k,:} =\mathbb{E}_{k}\left(\alpha_k\mC_{k}\right)  {=\sum_{k=1}^C \frac{1}{\alpha_k} \alpha_k \mC_k}\,, \qquad  {\text{with}\quad k \sim p_k=\frac{1}{\alpha_k},}
\end{equation}
where $\mC_k = \mA_{:,k}\mB_{k,:}$ is a rank-1 matrix formed from the outer product of $k$-th column of $\mA$ and the $k$-th row of $\mB$  {and the probabilities  $p_k$  can be chosen by the user.}

This rewriting expresses $\mC$ as the expectation of random rank-1 matrices. Consequently, it suggests a Monte Carlo (MC) approximation:
\begin{equation}\label{eqn:matrix_matrix_MC}
\mC = \mathbb{E}_k(\alpha_k \mC_{k}) \ \approx \ \frac{1}{c}\sum_{i=1}^c \alpha_{k_i}\mC_{k_i},
\end{equation}
where $\{k_i\}_{i=1}^c$ are i.i.d. samples drawn from the distribution $\{1/\alpha_j\}_j$. Equivalently, the approximation can be cast in linear algebra form:
\begin{equation}\label{eqn:matrix_matrix_MC_S}
\mC = \mA \mB \ \approx \ \mA  \mS^\top \mS \mB,
\end{equation}
where $\mS \in \R^{c \times C}$ is a sampling matrix with rows
\begin{equation}\label{eqn:matrix_matrix_MC_samplingMatrix}
\mS_{i,:} = \sqrt{\alpha_{k_i}}, {\bf e}_{k_i}^\top,
\qquad k_i \sim \{1/\alpha_j\}_j\,.
\end{equation}
See the sampling illustration in \cref{fig: hess sampling1} below.

\begin{center}
\resizebox{0.7\columnwidth}{!}
{\begin{tikzpicture}[
 line/.style ={draw, thick, -latex}                ]

\node[rectangle,minimum width=4.5cm, minimum height = 1.8cm, align=center,fill=orange,fill opacity=0.6,text opacity=1,font=\bfseries] at (-10,-1) (big block3) {} ;
\node[rectangle,minimum width=1.8cm, minimum height = 4.5cm, align=center,fill=orange,fill opacity=0.6,text opacity=1,font=\bfseries] at (-6,-1) (big block4) {} ;
\draw [decoration={brace,amplitude=4pt},decorate,line width=1pt] ($(big block4)+(2.5em,12ex)$) -- ($(big block4)+(2.5em,-14ex)$);
\node[rectangle,minimum width=2cm, minimum height = 2cm, align=center,fill=white,fill opacity=0,text opacity=1, font=\Large\bfseries] at (-4.5,-0.5) {$C$};
\node[rectangle,minimum width=2cm, minimum height = 2cm, align=center,fill=white,fill opacity=0,text opacity=1, font=\Large\bfseries] at (-10,-2.2) {$\bf A$};
\node[rectangle,minimum width=2cm, minimum height = 2cm, align=center,fill=white,fill opacity=0,text opacity=1, font=\Large\bfseries] at (-6,-3.5) {$\bf B$};

\node[rectangle,minimum width=3cm, minimum height = 2cm, align=center,fill=olive,fill opacity=0.5,text opacity=1,font=\bfseries] at (-0.5,-1) (block) {} ;
\node[rectangle,minimum width=2cm, minimum height = 3cm, align=center,fill=olive,fill opacity=0.5,text opacity=1,font=\bfseries] at (2.5,-1.5) (block1) {};
\draw [decoration={brace,amplitude=4pt},decorate,line width=1pt] ($(block1)+(3em,9ex)$) -- ($(block1)+(3em,-9ex)$);
\node[rectangle,minimum width=2cm, minimum height = 2cm, align=center,fill=white,fill opacity=0,text opacity=1, font=\Large\bfseries] at (4,-0.5) {$c$};
\node[rectangle,minimum width=2cm, minimum height = 2cm, align=center,fill=white,fill opacity=0,text opacity=1, font=\Large\bfseries] at (-0.5,-2.2) {$\bf A  S^\top $};
\node[rectangle,minimum width=2cm, minimum height = 2cm, align=center,fill=white,fill opacity=0,text opacity=1, font=\Large\bfseries] at (2.3,-3.3) {$\bf SB$};

\path[draw=black, line width=0.5mm, -{Triangle[length=4mm]}, shorten >=1mm, shorten <=0.5mm]   (-4.5,-1) to [left]    (-2.5,-1);
\end{tikzpicture}}
\vspace{-0.5em}

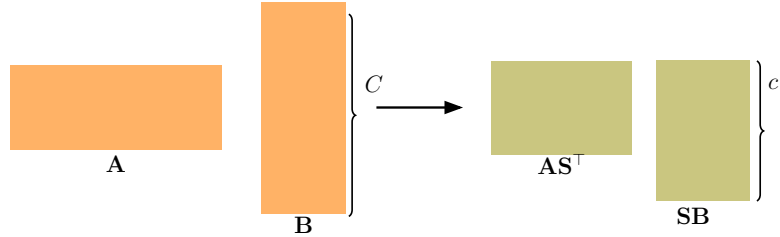
\captionof{figure}{Matrix-matrix product sampling demonstration.}
\label{fig: hess sampling1}
\end{center}

The quality and efficiency of this randomized solver depend on the choice of probabilities ${1/\alpha_k}$ and the number of samples $c$. Drineas et al.\cite{drineas2006fast} showed that the optimal configuration is achieved when
\begin{equation}\label{eqn:best_conf_mm_product} \alpha_k\propto |\mA_{:,k}||\mB_{k,:}| \,,\qquad\text{and}\qquad c=O\left( \frac{1+\log{\delta^{-1}}}{\epsilon^2}\right) 
\end{equation}
which guarantees that the Monte Carlo approximation\eqref{eqn:matrix_matrix_MC} satisfies
\begin{equation}\label{eqn:matrix_matrix_MC_accuracy}
|\mC - \mA \mS^\top \mS  \mB|_F < \epsilon
\qquad \text{with probability at least } 1-\delta.
\end{equation}
This randomized viewpoint is mathematically elegant: the matrix product is reconstructed as the average of randomly sampled rank-1 contributions. 
In practice the method  has computational benefits  only when the required number of samples $c$ is much smaller than $C$.

\subsubsection{Matrix-Matrix product for PDE based experimental design}\label{sec:MM_pde}
PDE-based experimental design provides a natural example of an effective application of the randomization approach for computing matrix-matrix product, through the  computation of the Fisher Information Matrix (FIM)~\eqref{eqn:hessian_L}.
Recall from~\eqref{eqn:hessian_L} that, in regions where $\vecsigma$ is close to the ground truth, the Hessian matrix of the loss functional is well approximated by the FIM. 
A distinctive feature of this matrix is that it decomposes as a sum of many rank-1 contributions:
\begin{equation}
\label{eq:Hessian}
    \mathrm{Hess}_{\vecsigma} L\approx {\bf H }= \mA^\top \mA  = \sum_{\theta} \mA_{\theta,:}^\top\mA_{\theta,:}\,.
\end{equation}
Here, each row $\mA_{\theta,:}$ corresponds to the Fr\'echet derivative $\frac{d}{d\vecsigma}\mathrm{pde}(\vecsigma,\theta)$ associated with a particular design $\theta$. In the semi-infinite setting with $\vecsigma \in \Sigma_n$ and $\theta\in\Theta$, this makes $\mA$ a matrix with $n$ columns and infinitely many rows, a scenario that fits seamlessly into the randomized matrix–matrix product framework of~\eqref{eqn:matrix_matrix_MC}.

Following the sampling formulation~\eqref{eqn:matrix_matrix_MC_S}-\eqref{eqn:matrix_matrix_MC_samplingMatrix}, one selects a subset of designs $\Theta_c \subset \Theta$ to construct a down-sampled Hessian:
\[
\mathbf{H}_S = \mA^\top \mS^\top \mS\mA\,,
\]
where $\mS$ encodes the randomized sampling. According to~\eqref{eqn:best_conf_mm_product}, the optimal strategy is to sample in proportion to $|\mA_{\theta,:}|^2$, i.e., the contribution of each design to the FIM.

A corollary of the analysis in~\cite{drineas2011faster}, as applied in~\cite{HKL24}, guarantees that with high probability:
\[
\lambda(\mathbf{H}_S)>\lambda(\mathbf{H}) -\epsilon \quad \text{with probability} \geq 1-\delta\,,
\]
where $\lambda$ denotes the smallest eigenvalue. This bound ensures the positivity of the sketched Hessian and thereby the well-posedness of the subsampled optimization problem.

\begin{figure}
\includegraphics[width = 0.24\textwidth]{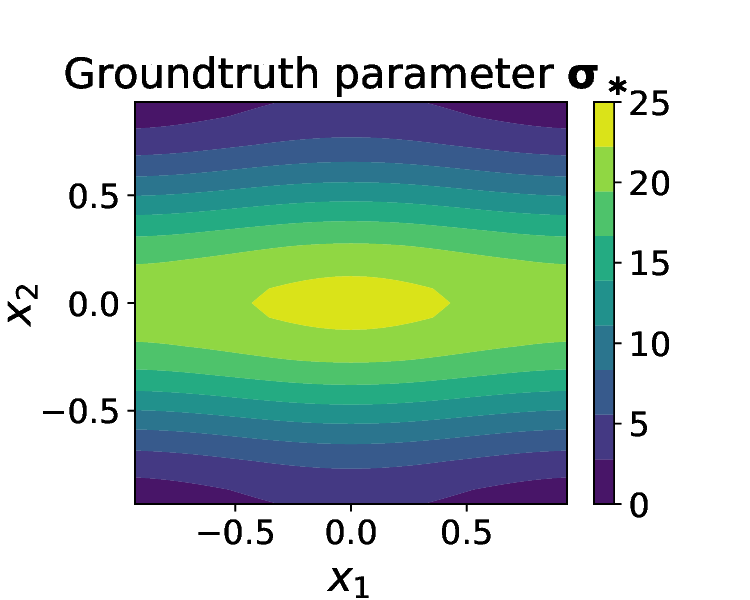}
\includegraphics[width = 0.24\textwidth]{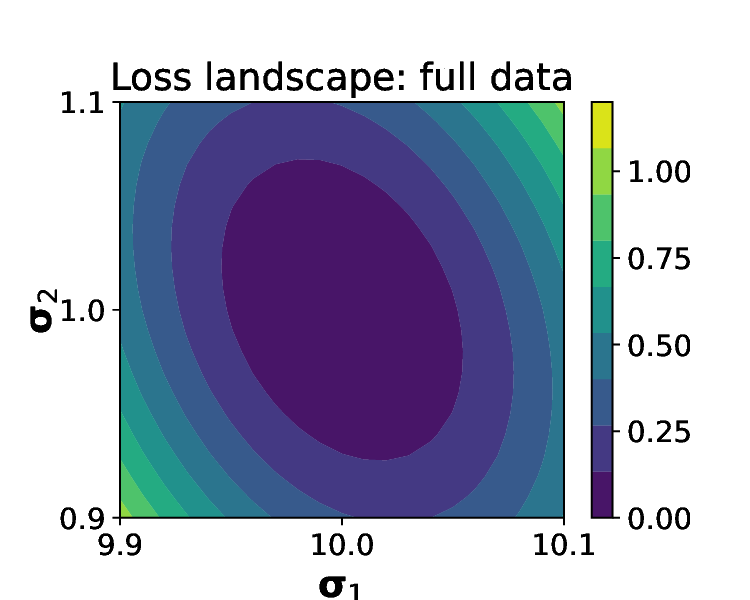}
\includegraphics[width = 0.24\textwidth]{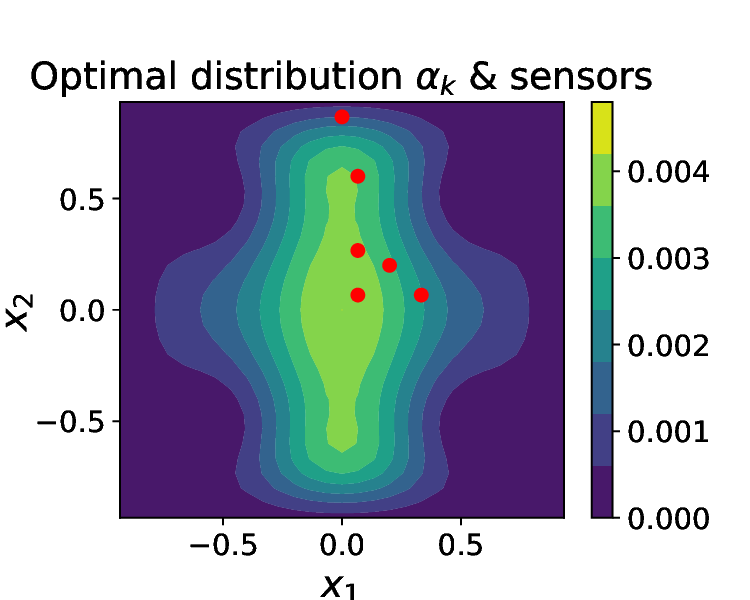}
\includegraphics[width = 0.24\textwidth]{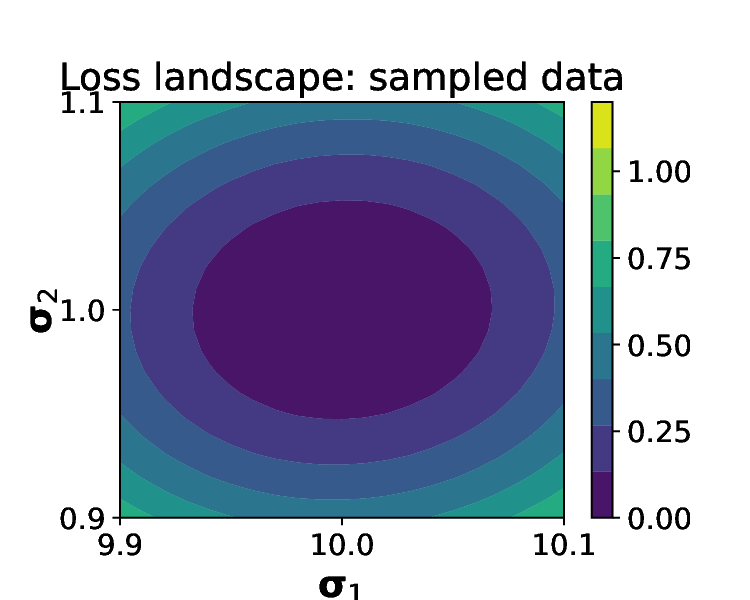}\\
\includegraphics[width = 0.24\textwidth]{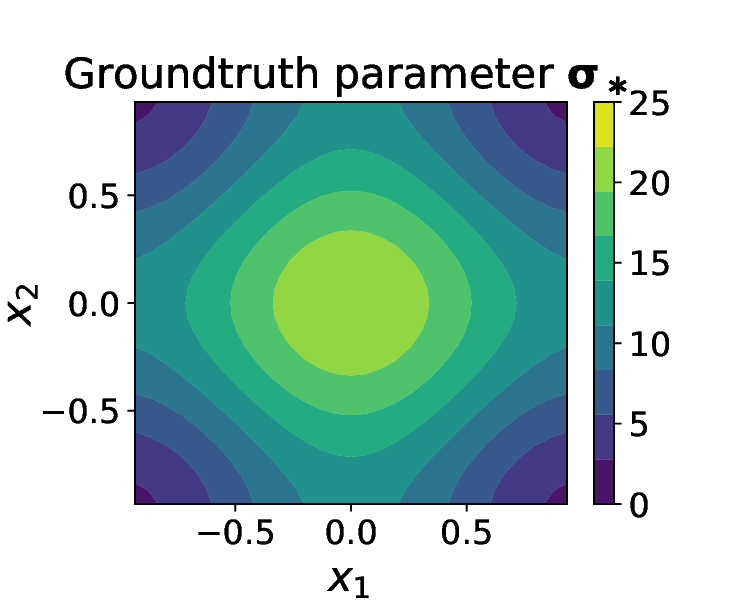}
\includegraphics[width = 0.24\textwidth]{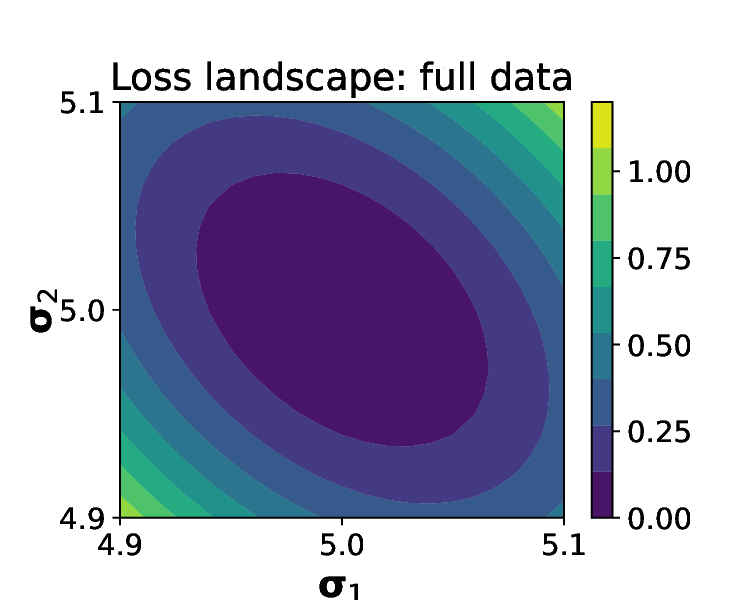}
\includegraphics[width = 0.24\textwidth]{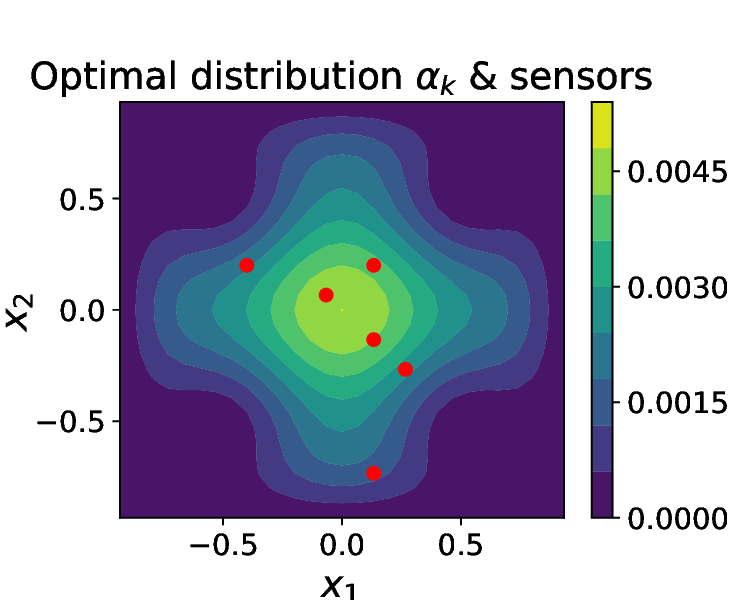}
\includegraphics[width = 0.24\textwidth]{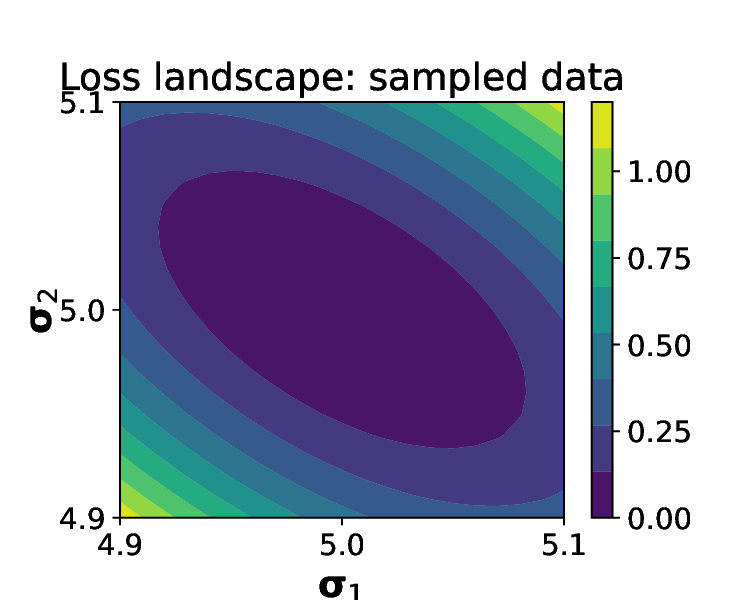}\\
\caption{Loss landscapes $L(\vecsigma)$ for two different ground truth parameters $\vecsigma_\ast$ (left column)  {with full data by setting with sensors placed everywhere in the domain $[-1,1]^2$} (middle left), and sampled data (right column)  {with $6$ sensors (red dots) randomly sampled} according to the optimal sampling distribution (middle right).  {The parameterization is $\vecsigma(x) = 12+\vecsigma_1 \cos(\pi x_1)+\vecsigma_2 \cos(\pi x_2)$.}}
\label{fig:LossLandscapes_Sampling}
\end{figure}

This randomized FIM approximation has been applied, for instance, to the linearized stationary Schr\"odinger potential reconstruction problem (\ref{itm: E3}).  {When fixing the boundary data $\phi_{\theta_1}$, the design space consists of all possible sensor locations $\theta_2\in[-1,1]^2$.} The results demonstrate that the  optimal  {distributions of} sensor locations are highly sensitive to the underlying ground truth $\vecsigma_\ast$, as illustrated in~\Cref{fig:LossLandscapes_Sampling}.  {The finite data setting that emerges from drawing a small number of sensor locations  from this distribution adheres good cost function convexity.}

This sensitivity is not a desired feature for reconstruction, and can potentially be overcome by a natural extension of this perspective to a sequential selection of experimental designs that involves feedback from previous experiments, as described in \Cref{sec:activeLearning}, which is left for future work. 
Moreover, a combination with a matrix completion approach in \Cref{sec: completion} promises enhanced reconstruction  efficiency  and may help bridge the gap between the semi-infinite and finite data regime.

\subsection{Column subset selection and conditioning guarantees}
Identifying a well-conditioned subset of columns from a large, wide matrix is another fundamental question in numerical linear algebra:  {Suppose we are given a ``dictionary" matrix $\mA\in\mathbb{R}^{c\times N}$, with $c\leq N$. Since $\mA$ only has $c$ rows, the row space is at most $c$ dimensional, meaning a small subset of columns in $\mA$ can potentially already contain all the $c$-dimensional information. Naturally, one would ask: } how can one select columns from $\mA$ to form a submatrix $\mA_\mS$ so that $\mA_\mS$ is well-conditioned? This problem is central to many applications in data science, including sparse signal and image reconstruction.

Mathematically, the selection process can be described using a sampling operator $\mS$, and the conditioning is evaluated via the Gram matrix
\begin{equation}\label{eqn:gram_column_selection}
{\bf H}_\mS \;=\; \mS^\top \mA^\top \mA \mS \;=\; \mA_\mS^\top \mA_\mS \;\in \R^{n \times n}\,, \quad \mbox{where $\mA_\mS:= \mA \mS$.}
\end{equation}
The goal is to design $\mS$ so that $\text{cond}({\bf H}_\mS) = \frac{\lambda_{\max}({\bf H}_\mS)}{\lambda_{\min}({\bf H}_\mS)}$ is small, where $\lambda_{\max,\min}$ denote the extreme eigenvalues. See illustration in Fig.~\ref{fig: hess sampling}.
\begin{center}
\resizebox{0.85\columnwidth}{!}
{\begin{tikzpicture}[
 line/.style ={draw, thick, -latex}                ]

\node[rectangle,minimum width=4cm, minimum height = 4cm, align=center,fill=orange,fill opacity=0.6,text opacity=1,font=\bfseries] at (-16,0) (big block1) {$\bf H$} ;
\node[rectangle,minimum width=2cm, minimum height = 2cm, align=center,fill=white,fill opacity=0,text opacity=1, font=\Large\bfseries] at (-16,-2.5) {ill-conditioned};
\node[rectangle,minimum width=2cm, minimum height = 4cm, align=center,fill=white,fill opacity=0,text opacity=1,font=\bfseries] at (-13,0) {\large{=}};
\node[rectangle,minimum width=3cm, minimum height = 4cm, align=center,fill=orange,fill opacity=0.6,text opacity=1,font=\bfseries] at (-11,0) (big block1) {${\bf A}^\top$} ;
\node[rectangle,minimum width=4cm, minimum height = 3cm, align=center,fill=orange,fill opacity=0.6,text opacity=1,font=\bfseries] at (-7,0.5) (big block2) {$\bf A$} ;

\node[rectangle,minimum width=0cm, minimum height = 0cm, align=center,fill=white,fill opacity=0,text opacity=1,font=\Large\bfseries] at (-4.4,0.5) {$c$};
\node[rectangle,minimum width=0cm, minimum height = 0cm, align=center,fill=white,fill opacity=0,text opacity=1,font=\Large\bfseries] at (-7,2.5) {$N$};

\node[rectangle,minimum width=3cm, minimum height = 2cm, align=center,fill=olive,fill opacity=0.5,text opacity=1,font=\bfseries] at (0,0) (block) {$\bf A_S^\top$} ;

\draw [decoration={brace,amplitude=4pt},decorate,line width=1pt] ($(big block2)+(6em,8ex)$) -- ($(big block2)+(6em,-10ex)$);
\draw [decoration={brace,amplitude=4pt},decorate,line width=1pt] ($(big block2)+(-5em,10ex)$) -- ($(big block2)+(5em,10ex)$);
\node[rectangle,minimum width=2cm, minimum height = 3cm, align=center,fill=olive,fill opacity=0.5,text opacity=1,font=\bfseries] at (3,0.5)(big block3) {$\bf A_S$};
\node[rectangle,minimum width=2cm, minimum height = 4cm, align=center,fill=white,fill opacity=0,text opacity=1,font=\bfseries] at (4.5,0) {\large{=}};
\node[rectangle,minimum width=2cm, minimum height = 2cm, align=center,fill=gray,fill opacity=0.5,text opacity=1,font=\bfseries] at (6,0)  {$\bf H_S$};

\path[draw=black, line width=0.5mm, -{Triangle[length=4mm, bend]}, shorten >=1mm, shorten <=0.5mm]    
        (-4.5,0) to    (-2,0);
\node[rectangle,minimum width=2cm, minimum height = 2cm, align=center,fill=white,fill opacity=0,text opacity=1, font=\Large\bfseries] at (6,-2) {well-conditioned};
\node[rectangle,minimum width=0cm, minimum height = 0cm, align=center,fill=white,fill opacity=0,text opacity=1,font=\Large\bfseries] at (3,2.5) {$n$};
\draw [decoration={brace,amplitude=4pt},decorate,line width=1pt] ($(big block3)+(-2.5em,10ex)$) -- ($(big block3)+(2.5em,10ex)$);
\node[rectangle,minimum width=0cm, minimum height = 0cm, align=center,fill=white,fill opacity=0,text opacity=1,font=\Large\bfseries] at (4.5,0.5) {$c$};
\draw [decoration={brace,amplitude=4pt},decorate,line width=1pt] ($(big block3)+(3em,8ex)$) -- ($(big block3)+(3em,-10ex)$);
\end{tikzpicture}}
\vspace{-0.5em}
\captionof{figure}{Hessian (column) sampling: ${\bf H} = \mA^\top \mA \in \R^{N \times N} \rightarrow {\bf H}_\mS = \mA_\mS^\top \mA_\mS \in \R^{n \times n}.$}
\label{fig: hess sampling}
\end{center}

Randomization provides powerful tools for this problem. A key result of~\cite{tropp2008conditioning} gives probabilistic guarantees on conditioning when columns are sampled uniformly. Assuming the columns of $\mA$ are normalized, the number of columns $n$ that can be selected while retaining good conditioning satisfies
\begin{equation}
\label{eqn: column sample1}
n \;\leq\; \min\!\left( o\!\left(\frac{1}{\mu^2 \log c}\right), \;\; o\!\left(\frac{N}{\| {\bf H}\|_2} \right) \right),
\end{equation}
with high probability, where
\begin{equation}
\label{eqn: condition concentration1}
\mathbb{P}\!\left( \text{cond} ({\bf H}_\mS) \;\leq\; \frac{1 + \tau(n)}{1 - \tau(n)} \right) \;\geq\; 1 - \frac{1}{c}.
\end{equation}
Here, $\mu = \max_{i \neq j} \vert \langle \mA_{:, i}, \mA_{:, j} \rangle \vert$ is the coherence between columns, and $\tau(n)$ is a distortion parameter scaling as $\tau(n) = \mathcal{O}\!\left( \sqrt{\mu^2 n \log c} \;+\; n\frac{\|{\bf H} \|_2}{N}  \right)$.

The coherence $\mu$ plays a decisive role: If $\mu$ is large, columns are highly overlapping, and only a relatively small number of them can be selected while preserving conditioning.  Therefore, smaller $\mu$ enlarges the allowable range of $n$, consistent with the intuition that independence between columns enables more extensive sampling without loss of stability.  {We discuss the second term $o\left(\tfrac{N}{\|{\bf H}\|_2}\right)$ of \eqref{eqn: column sample1} a little further. 
The inverse of this bound  controls the discrepancy between the sampled Hessian $\bf H_S$ and the identity $\bf I$, specifically, $\mathbb{E} \| {\bf H_S} - {\bf I}\|_2 \leq n\frac{\| {\bf H}\|_2}{N}$. 
This bound in turn controls the condition number estimation in \eqref{eqn: condition concentration1}. 
Further, since $N = \mathrm{Tr}({\bf H}) = \sum_i \lambda_i({\bf H}) \;\leq\; c \,\|{\bf H}\|_2$, the ratio $\tfrac{N}{\|{\bf H}\|_2}$ is bounded by $c$, so the second term of \eqref{eqn: column sample1} implies that $n \le o(c)$, which in itself is not a particularly useful bound.}
In summary, the sampling dimension $n$ \eqref{eqn: column sample1} scales  {sublinearly} with the matrix rank $c$ and is influenced by the factor of $1/\mu^2.$

\subsubsection{Column subset selection for PDE based experimental design}
As discussed in~\eqref{eqn:hessian_L} and Section~\ref{sec: sampling}, the Hessian of the objective near the ground truth is dominated by the FIM. Since the FIM has an outer-product structure closely resembling~\eqref{eqn:gram_column_selection}, it is natural to expect that column subset selection strategies can provide useful guidance for data selection in PDE-based inverse problems.

From~\eqref{eq:Hessian}, the Hessian is expressed as an outer product of $\mA$, with rows $\mA_{\theta,:}=\frac{\dd}{\dd \vecsigma}\mathrm{pde}(\theta,\vecsigma) \in \R^N$. Applying the column selection framework in~\eqref{eqn:gram_column_selection} corresponds to fixing a set of $\theta$ samples and restricting $\vecsigma$. This is the opposite setting from Section~\ref{sec:MM_pde}: Here, the experiments are already collected, and the task is to identify subsets of parameters $\vecsigma$ that can be reconstructed in a well-conditioned manner. 
This strategy corresponds to \emph{random sampling in parameter space}.

Direct application of prior results~\eqref{eqn: column sample1}–\eqref{eqn: condition concentration1} is obstructed by the fact that the columns of $\mA$ are not normalized in our application. 
Nevertheless, by adapting proof techniques from~\cite{RV08,tropp2008conditioning}, the work~\cite{JLNS24} generalizes the guarantees to PDE settings and establishes that
\begin{equation}
\label{eqn: column sample}
n \leq \min \left( o\left(\frac{\text{Tr} ({\bf H})}{\| {\bf H}\|_2} \right), o\left(\frac{1}{\mu^2\, \log c} \frac{\text{Tr}({\bf H})^2}{\| {\bf H}\|_F^2} \right) \right),
\end{equation}
with the probabilistic conditioning bound
\begin{equation}\label{eqn: condition concentration}
\mathbb{P} \left(\text{cond} ({\bf H_S}) \leq \frac{L + \tau(n)}{\ell - \tau(n) }\right)  \geq 1 - \frac{1}{c}\,.
\end{equation}
Here the distortion quantity is given by $\tau(n) = e^{\frac{1}{4}} \left(2n \frac{\|{\bf H} \|_2}{\text{Tr} ({\bf H})} + 12 \mu \sqrt{n  \log c}\, \frac{\| {\bf H}\|_F}{\text{Tr}({\bf H})}\right)$, and the normalization constants account for extreme column norms:
\[
\ell: = \frac{N}{\text{Tr} ({\bf H})} \min_{i=1\dots N} \| \mA_{:,i}\|_2^2, \quad L: = \frac{N}{\text{Tr} ({\bf H})} \max_{i=1\dots N} \| \mA_{:,i}\|_2^2\,.
\]
The coherence parameter must also be adjusted to account for non-uniform column magnitudes: $\mu: = \frac{N}{\| {\bf H}\|_F} \max_{i \neq j = 1\dots N} \vert \langle \mA_{:,i}, \mA_{:,j} \rangle \vert$.

The second term in~\eqref{eqn: column sample} highlights the role of coherence: if two parameter directions $\vecsigma_1,\vecsigma_2$ yield nearly parallel sensitivity vectors $\tfrac{\dd}{\dd \vecsigma}\mathrm{pde}(:,\vecsigma_1)$ and $\tfrac{\dd}{\dd \vecsigma}\mathrm{pde}(:,\vecsigma_2)$, the effective number of reliably reconstructible parameters decreases.

\begin{figure}[!htb]
\centering
\subfloat[Visualization of $\mA$ and ${\bf H}$ from the EIT model.]{\includegraphics[width=0.5\linewidth]{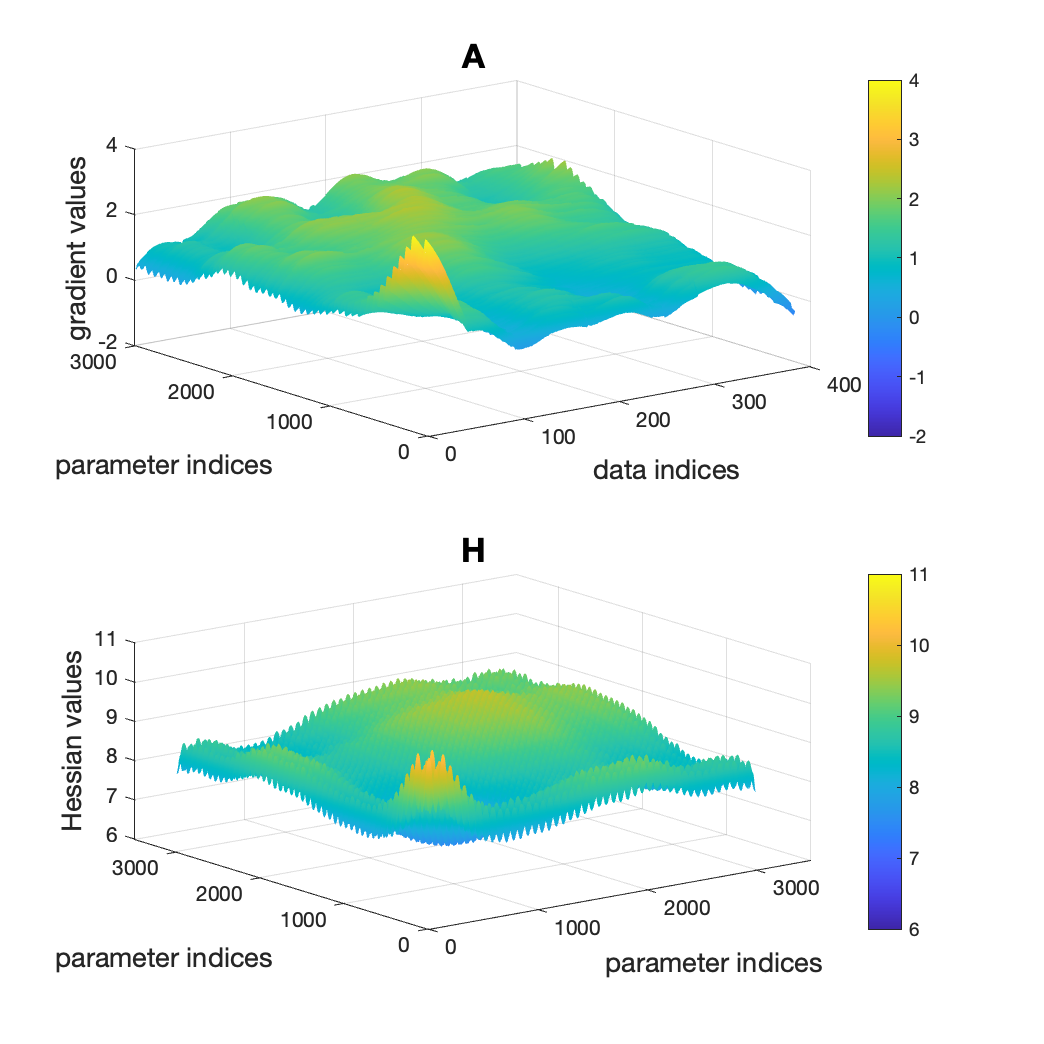}}~~
\subfloat[Sample complexity $n$ (horizontal) vs.\ condition number (vertical). The shaded region shows the 20\%–80\% quantile range over 10,000 simulations.]{\includegraphics[width=0.48\linewidth]{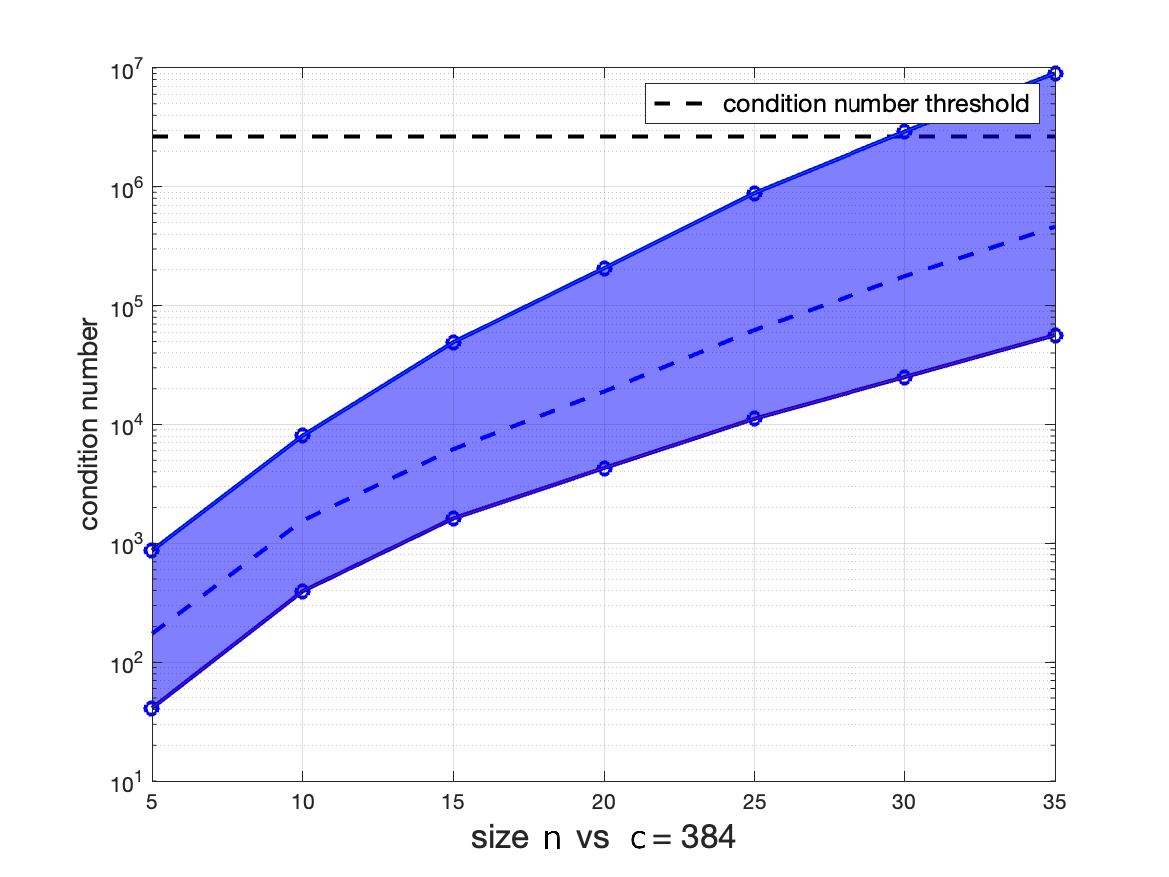}}
\caption{Column subset selection in PDE-based inverse problems.}
\label{fig:hess_sampling_pde}
\end{figure}

Numerical experiments validate these results. In Fig.\ref{fig:hess_sampling_pde}, $\mA$ is generated from the linearized EIT model~\ref{itm: E2} with $c=384$. 
As $n$ increases, the condition number of ${\bf H_\mS}$ grows consistently with the concentration inequality~\eqref{eqn: condition concentration}. For reference, the baseline threshold $\tfrac{L}{\ell}$ is also shown (dashed line).

To conclude, understanding unique reconstruction and well-posedness from limited data information is a key challenge in PDE inverse problems. 
Random sampling in the parameter space can be used to build a probabilistic framework and discern the number of parameters identifiable in a well-conditioned manner. 
Open research avenues for future study include: 1. Incorporating greedy algorithm and importance sampling strategy to maximize the number of parameters recovered; 2. Extending the analysis to the infinite-dimensional PDE setting with full parameter space $\Sigma = \lim_{N \to \infty}\Sigma_N$, motivated by the concentration results \eqref{eqn: column sample}-\eqref{eqn: condition concentration} being independent from the ambient parameter dimension $N$.

\subsection{Randomized SVD}
In singular value decomposition (SVD), another classical task in numerical linear algebra, we start with a rectangular matrix $\mathbf{C}\in\R^{m\times n}$ of rank $d \le \min(m,n)$, and find matrices $\mathbf{U}$, $\mathbf{S}$, $\mathbf{V}$ such that
\[
\mathbf{C}=\mathbf{U}\mathbf{S}\mathbf{V}^\top = \sum_is_i\mathbf{U}_{:,i}\left(\mathbf{V}_{:,i}\right)^\top\,,
\]
where $\mathbf{S}=\text{diag}\{s_i\}$ collects singular values that is ordered in a descending manner, while $\mathbf{U}\in\R^{m\times d}$ and $\mathbf{V}\in\R^{n\times d}$ have orthonormal columns, which are the left and right singular vectors or $\mathbf{C}$, respectively {.}
Classically, computation of the  SVD  starts with a bidiagonalization via Householder transformations, which transforms $\mathbf{C}$ into a matrix that has two diagonal components. 
The total cost is $O(\max(m,n)\cdot\min(m,n)^2)$.

Randomized SVD (RSVD) returns estimates to $\mathbf{U}$, $\mathbf{V}$ and $\{s_i\}$ with high accuracy and probability. 
As thoroughly analyzed in the foundational works \cite{MARTINSSON201147, liberty2007randomized, halko2011finding}, it is a much more cost-efficient algorithm than classical SVD for large, dense matrices, and is highly parallelizable and memory-efficient. 
The total cost is about $O(mnd)$, and in comparison with classical solvers, the saving is most pronounced when the large matrix $\mathbf{C}$ is known to be of low rank, that is, $d\ll\min\{m,n\}$. 
As a consequence, RSVD finds applications in domains that traditionally rely on Principal Component Analysis (PCA) \cite{rokhlin2010randomized}, such as bioinformatics \cite{darnell2017adaptive} and image processing \cite{https://doi.org/10.1155/2012/409357}. 

RSVD has been successfully deployed for experimental design purposes for PDE-based inverse problems as well, especially through the optimal design angle. 
Recalling~\eqref{eqn:optimal_design}, the task at hand it to find 
 {the design $\Theta_c\subset \Theta$} 
so that the induced variance matrix  {$\hat \Gamma$} 
is small, either in terms of largest eigenvalue ($\Phi_{E}$), or the summation of all eigenvalues ($\Phi_A$) or the multiplication ($\Phi_D$), all of which require computation of eigenvalues. 
Considering the size of  {$\hat \Gamma$}, RSVD that returns all eigenvalues efficiently becomes very useful. 
This perspective was adopted in~\cite{SAI17,aarset2025global, herman2020randomization, A21,wu2023fast,APSG14}.

\subsection{Compressed sensing}
\label{sec: CS}
Compressed sensing (CS) is a signal processing technique that leverages randomized numerical linear algebra ideas, and it has promising applications in PDE-based inverse problems. 
As its name suggests, CS acquires information about a signal in a compressed manner. 
When the object to be reconstructed lies in a high-dimensional ambient space $\R^N$ but is known to be sparse (that is, supported on only a small fraction of the coordinates), the sparsity can be exploited to drastically reduce the number of measurements required. 
In particular, if the sensing mechanism satisfies certain structural properties, the number of ``sensors'' scales with the intrinsic information content $n$ rather than with the ambient dimension $N$.

A standard reconstruction method takes the form of an optimization problem:
\[
\min\|\mathbf{x}\|_1\,,\quad\text{such that}\quad \mathbf{P}_{\Omega}\mathbf{U}\cdot\mathbf{x} = \mathbf{P}_\Omega\mathbf{y}\,.
\]
Here $\mathbf{U}$ is an orthonormal sensing matrix collecting all possible measurements, and $\mathbf{P}_\Omega$ is a projection operator selecting entries according to the mask $\Omega$ of size $|\Omega|=c$. 
Because $\mathbf{x}$ is assumed sparse, the $\ell_1$ norm promotes sparsity in the recovered solution. 
Although the system is highly underdetermined, posing $c \ll N$ constraints for $N$-entry of reconstruction, exact recovery is still possible when $\mathbf{U}$ is incoherent with the coordinate basis. 
Classical results show that if $c \gtrsim n\log N$, then exact recovery of $\mathbf{x}$ is achievable with high probability~\cite{candes_cs,donoho_cs,Candes_2007_fourier}.

These foundational results, originally developed in finite-dimensional linear algebra, extend naturally to infinite-dimensional settings that arise in PDE-based inverse problems. In this context, the to-be-reconstructed object is a function $f(x)$ rather than a finite-dimensional vector, and the sensing matrix becomes a continuous operator $\mathcal{U}$. Nevertheless, if $f(x)$ is sparse in a basis system $\mathcal{D}=\{\phi_i\}$ and if $\mathcal{D}$ and $\mathcal{U}$ are mutually incoherent, compressed sensing theory ensures stable recovery even from finitely many measurements~\cite{ALBERTI2021105,Hansen}.

\subsection{Matrix completion}
\label{sec: completion}
Matrix completion seeks to complete a dense matrix $\mA\in\R^{C\times C}$, given certain entries $\mA_{ij}$ for $(i,j)\in\Omega$, a sparse mask. 
The task is applicable to recommendation system~\cite{KBV09} and genomics~\cite{CZ16}. The problem itself does not necessarily belong to RNLA, but some solvers use random matrices and the associated techniques. 
In its generic form, there is no hope  to reconstruct a matrix with a few entries, but if certain assumptions are satisfied, the reconstruction problem can be solved. Notably, when the matrix is known a-priori to be of low rank  and its column and row space known to be decoherent, the following optimization problem returns an optimizer that is known to recover $\mA$ with high probability\cite{CR08,CT10,R11}.
\begin{equation}
\label{eqn: nuclear norm min}
\begin{array}{lll}
\text{minimize} &\| \mX \|_{\ast} & \\
\text{subject~to~}& \mX_{ij} = \mA_{ij}, & \forall~(i,j) \in \Omega\,.
\end{array}
\end{equation}
where $\Omega$ is a randomly sampled mask with each entry following the Bernoulli distribution at a low sampling rate, and $\|\cdot\|_*$ is the nuclear norm (the sum of all singular values). 
Intuitively, the objective function is the $L_1$ norm version of singular values, in analogy to compressed sensing, where in both context it serves as a proxy for the $L_0$ norm. 
The minimization of $L_0$ norm means one looks for a reconstruction that has smallest rank and agrees with the given data. 
A crude conclusion suggests a sampling complexity of  $c\approx Cr$ where $r$ is the approximate rank.

It is easy to imagine these completion solvers being useful for qualitative experimental design associated with inverse problems. 
An inverse problem is to use ItO data pairs to reconstruct unknown parameters. 
Restricted by sampling and experiment budget, one does not have the full set of input-to-output data pairs, but rather a subset of entries. 
Can one use this subset to reconstruct the full data set, which is then feed into the true reconstruction?

\begin{figure}[!htb]
    \centering   \includegraphics[width=0.7\linewidth]{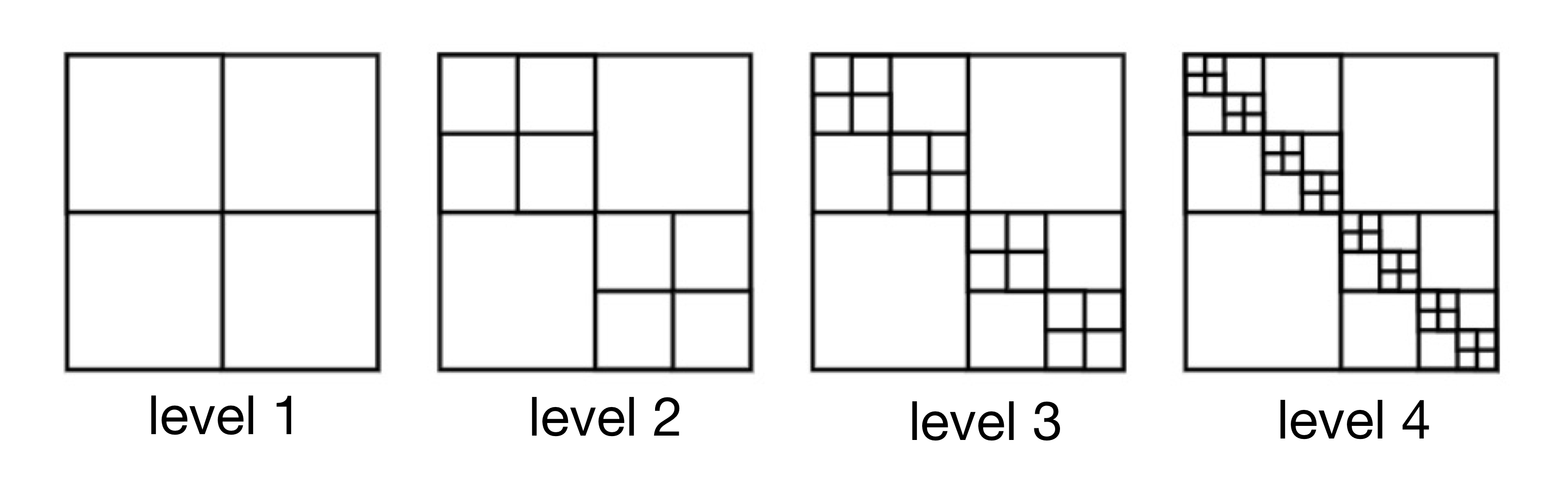} \includegraphics[width=0.16\linewidth]{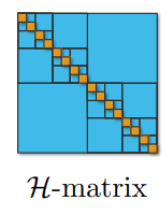}
    \vspace{-1em}
    \caption{Hierarchical matrix block structure on varying partition levels  and rank structure of the $\mathcal H$-matrix: the diagonal full-rank blocks (orange)  vs. blue low-rank off diagonal blocks.}
    \label{fig: hierachical structure}
\end{figure}

The approach is taken in~\cite{BLZ22} for the EIT problem, \ref{itm: E2} from section~\ref{sec:representatives}, for which the ItO map is the Dirichlet-to-Neumann data. It is important to note that a well-posed inverse problem typically corresponds to a full-rank data matrix that stores the full ItO map, which prevents immediate application of matrix completion method~\eqref{eqn: nuclear norm min}. 
Indeed, in~\cite{BLZ22}, recognizing the off-diagonally low rank feature of the DtN map, a dyadic decomposition on the diagonal blocks can be used to peel off the off-diagonal blocks as the partition level increases \cite{LLY11}, yielding a block structure as displayed in illustration \cref{fig: hierachical structure}.
Then matrix completion~\eqref{eqn: nuclear norm min} is applied to each off-diagonal block. 
The recovered full dataset is then deployed in PDE-constrained optimization for the execution of the inverse problem.

\section{Outlook}\label{sec:outlook}
This review has concentrated on the role of randomized numerical linear algebra (RNLA) in data selection for PDE-based inverse problems. While powerful and increasingly influential, our focus is necessarily narrow. 
At its core, RNLA is a linear framework whereas data selection is intrinsically nonlinear. 
As such, relying exclusively on RNLA overlooks essential aspects of the problem. 
Even within the linear regime, alternative toolboxes from modern machine learning and applied mathematics offer complementary strategies that may also prove valuable. 
In what follows, we highlight two such directions and discuss how they may broaden the scope of data selection methodologies.

\subsection{Active learning}\label{sec:activeLearning}

Most techniques described above are effectively \emph{passive}: The algorithm is prescribed once, and the user simply waits for the output.  {In contrast to that, real world data collection often happens sequentially, and feedback from previous experiments can be valuable in for future decisions.}

This observation motivates \emph{interactive} or \emph{adaptive} data selection: Begin with an initial dataset, perform a reconstruction, and then collect additional data informed by the current state of knowledge, repeating the process as needed.
\begin{figure}[!htb]
    \centering
    \begin{tikzpicture}[
    every node/.style={font=\sffamily},
    box/.style={
        draw,
        rounded corners,
        minimum width=2.2cm,
        minimum height=1cm,
        align=center,
        thick
    },
    arrow/.style={-{Latex[length=3mm]}, thick}
]

\def\radiusint{1.5cm} 
\def\radius{1.7cm}    
\def\offset{40} 

\node[box, fill=blue!20, draw=blue!60!black]       (model)       at (90:\radiusint)  {Model};
\node[box, fill=orange!25, draw=orange!80!black]   (design)      at (330:\radiusint) {Design};
\node[box, fill=green!25, draw=green!70!black]     (measurement) at (210:\radiusint) {Measurement};

\draw[arrow]
  ({90-\offset}:\radius)
  arc[start angle={90-\offset}, delta angle={-120+2*\offset}, radius=\radius]
  node[midway, above right]{\textbf{OED}};

\draw[arrow]
  ({330-\offset}:\radius)
  arc[start angle={330-\offset}, delta angle={-120+2*\offset}, radius=\radius]
  node[midway, below]{\textbf{Run Experiment}};

\draw[arrow]
  ({210-\offset}:\radius)
  arc[start angle={210-\offset}, delta angle={-120+2*\offset}, radius=\radius]
  node[midway, left]{\textbf{Parameter Inference}};

\end{tikzpicture}
    \caption{Greedy sequential optimal experimental design.}
    \label{fig:sOED}
\end{figure}
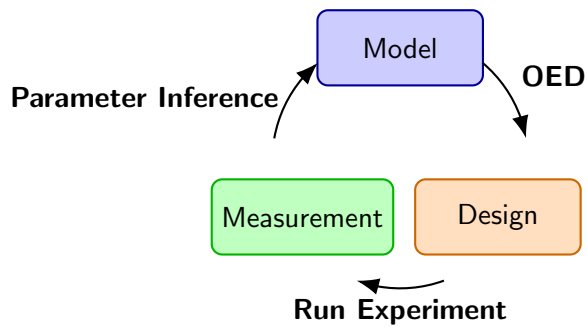

 {This paradigm allows one to select data sequentially in the most effective way. It not only reduces computational and sample complexity of data selection process, but also has the potential to improve the reconstruction \cite{elata2024adaptive}.}

The primary limitation of RNLA-based approaches is that they operate strictly within linear spaces, while experimental design is fundamentally nonlinear. 
As discussed in Section~\ref{sec:linearization}, the information value of an experiment depends strongly on the underlying medium $\vecsigma^*$. 
This dependence is evident in Figure~\ref{fig:LossLandscapes_Sampling}, where the weights assigned to measurements vary dramatically with different $\vecsigma^*$. 
A natural measure of data value is the sensitivity of the reconstruction to the data, often encoded in the Jacobian. 
In the RNLA framework, the Jacobian is treated as a fixed matrix.
In the nonlinear regime, however, the Jacobian itself evolves as reconstruction improves, shifting the sensitivity landscape and altering the relative value of new data.  {Adaptive approaches are more suitable to accommodate this behavior.}

Several concrete strategies capture this philosophy. 
One is to pose the problem as a bilevel optimization: The inner layer solves the reconstruction problem for a given dataset, while the outer layer selects the dataset itself~\cite{RCC18, RJ}. 
Another approach is greedy data acquisition~\cite{CABN19, W70, jagalur2021batch, wu2023fast}, in which data are acquired sequentially, each step exploiting information from previous ones, as described by \Cref{fig:sOED}. 
More sophisticated  approaches additionally allow for coordination of subsequent experiments and pose the problem in the reinforcement learning regime \cite{huan2016sequential,santosa2022bayesian}.

Crucially, this challenge is not unique to PDE-based inverse problems. Similar difficulties arise across machine learning tasks constrained by costly data acquisition, where \emph{active learning} has emerged as a central paradigm. We expect that cross-fertilization between active learning and nonlinear experimental design for PDEs will play an increasingly important role in the coming years~\cite{SCMW23}.

\subsection{Operating in infinite-dimensional settings}
Even within the linear regime, RNLA methods ultimately reduce large {\em but finite} linear systems via sampling or sketching. The design space is therefore still finite, no matter how large it is taken. 
This discreteness is still incompatible with PDE-based inverse problems, where the design space is continuously indexed and thus infinite-dimensional.

\begin{figure}[!htb]
\centering
\includegraphics[scale=0.42]{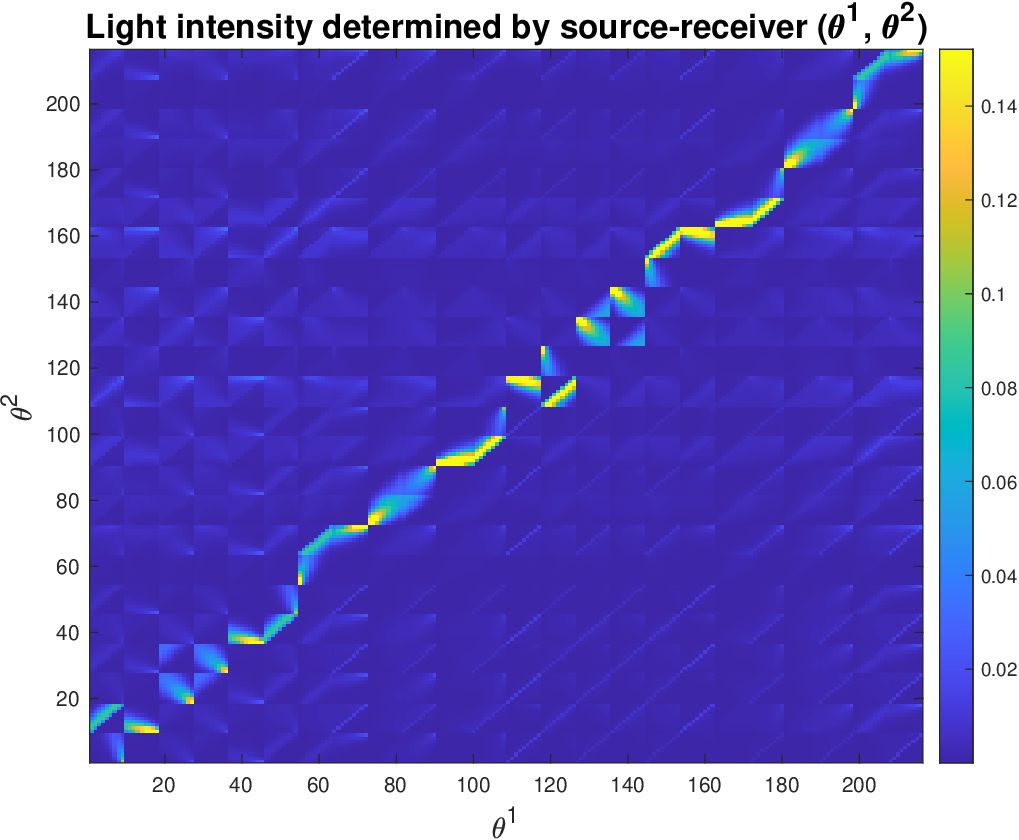}
\caption{Simulation of the radiative transfer equation characterizing laser beam propagation. The column index $\theta_1$ denotes the laser location and velocity direction. For each $\theta_1$, the rows correspond to boundary measurements at detector locations $\theta_2$. Only a few values of $\theta_2$ produce nontrivial readings, forming a measure-zero set in the full design space.}
\label{fig:LV}
\end{figure}

The importance of continuous indexing is illustrated in Fig.~\ref{fig:LV}, which simulates radiative transfer~\ref{itm: E4}. 
For each experiment (each column), only a handful of detector locations record nontrivial light intensity. 
The relevant information thus arises in a measure-zero subset of the design space. 
Any fixed discretization, no matter how fine, risks missing these measure-zero structures and therefore producing a poor design space for data selection.

Analogous phenomena arise in  {the Darcy flow variation of example~\ref{itm: E2}} and the Lorenz-63 problem~\ref{itm: E1}.
\begin{figure}[!htb]
\centering
\includegraphics[scale=0.38]{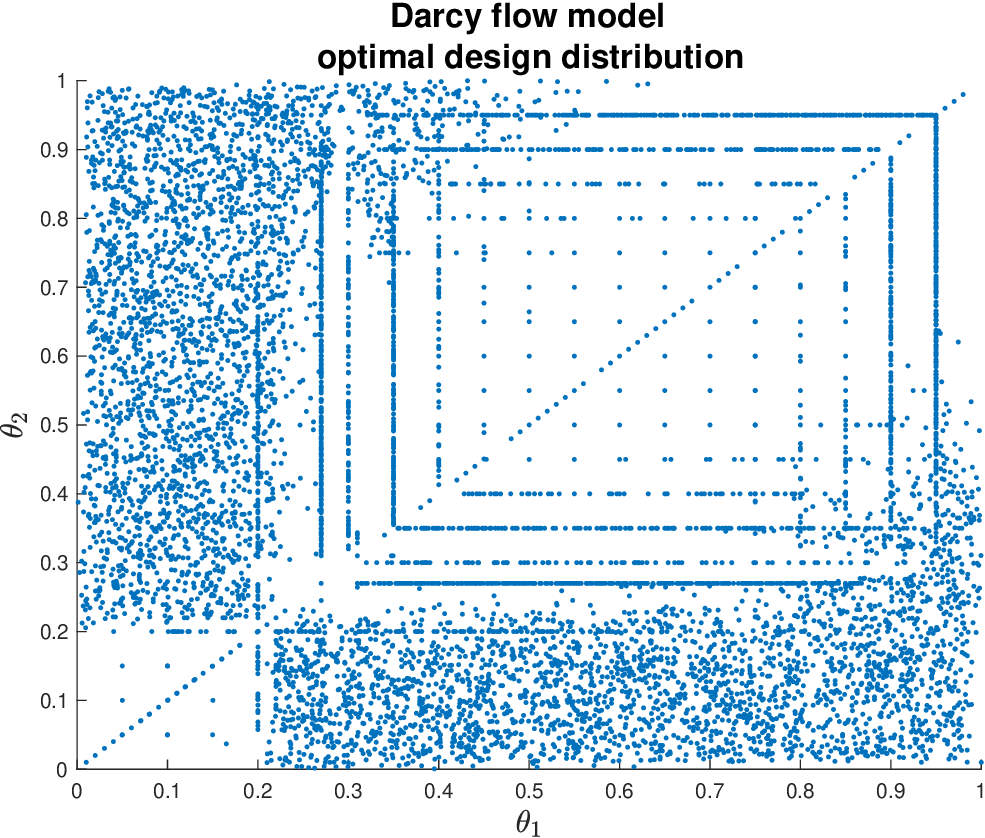}~~
\includegraphics[scale=0.4]{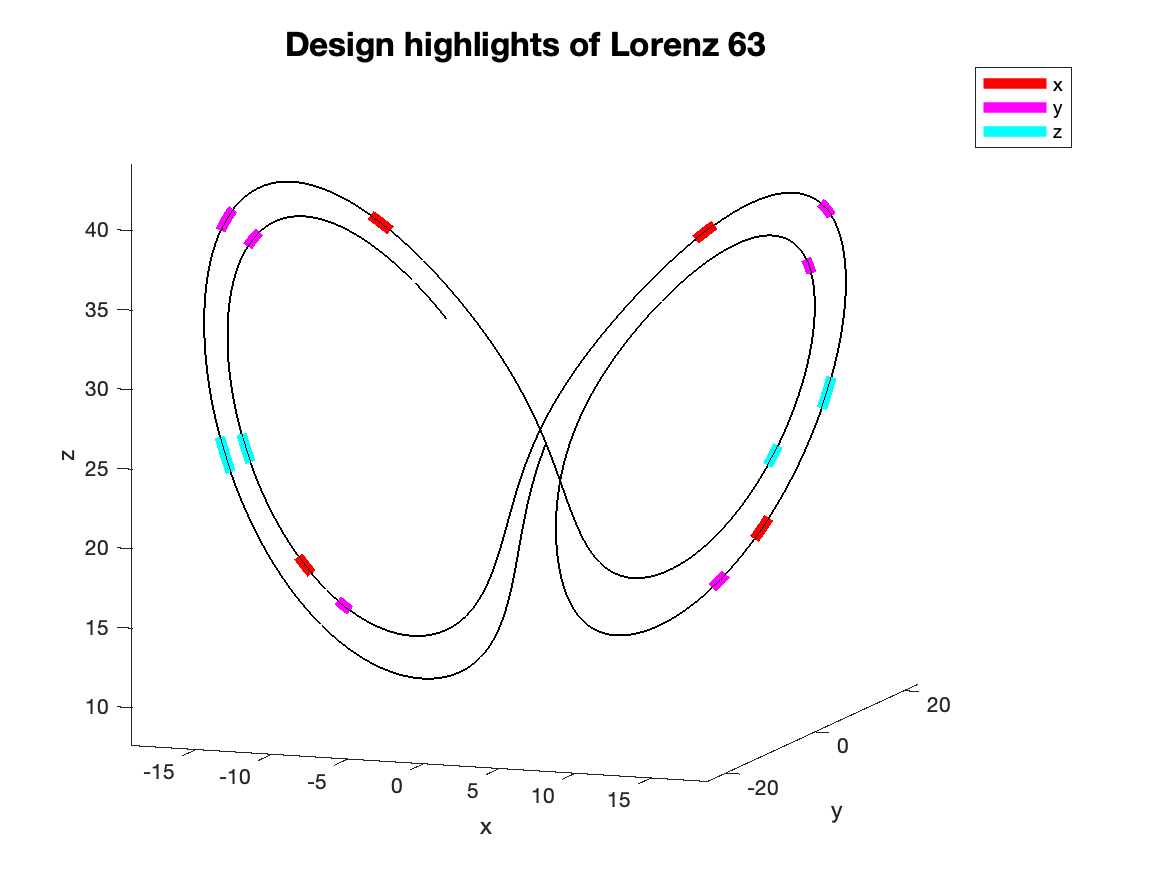}
\caption{Left panel: Optimal distribution over the continuous design space for the Darcy flow problem, showing extreme concentration of the data importance. Right panel: Important observation times in colored slots for parameter inference in the Lorenz-63 model. Different colors indicate different measured components $(x,y,z)$.}
\end{figure}

In such cases, methods that \emph{directly} operate in infinite-dimensional spaces become essential. 
Traditionally, such an approach was regarded as impossible, since most optimization solvers operate in $\R^d$ and exhibit complexity growing at least algebraically with $d$, rendering them infeasible in the infinite-dimensional limit. 
Recent advances have changed this landscape. 
Tools from optimal transport and gradient flows~\cite{AGS05, FG21, computationlOT} have clarified the structure of optimization over probability measure spaces, which are intrinsically infinite-dimensional. 
Techniques such as the mean-field limit and propagation of chaos~\cite{Louis-Pierre, fournier2015rate} further enable practical computation by translating infinite-dimensional PDEs into coupled ODE systems. 
These ideas, widely used in machine learning~\cite{Chizat, nitanda2025propagation, Nishikawa_2023, geshkovski2023the, chen2025transformer}, have only recently been adapted to experimental design problems~\cite{JGLW24}.

One such approach replaces discrete weights $w(\theta)$ in the classical optimal design problem~\eqref{eqn:optimal_design_weight} with a probability measure $\rho(\theta)$:
\begin{equation}
\begin{aligned}
\text{A-optimal}: &~ \Phi_A[\rho] = \mathrm{tr}\!\left( \left( \mA^\top\mA[\rho] \right)^{-1}\right), \\
\text{D-optimal}: &~ \Phi_D[\rho] = \det\!\left( \mA^\top\mA[\rho]\right),
\end{aligned}
\label{eqn:oed_criteria2}
\end{equation}
where
\[
\mA^\top\mA[\rho] = \int_{\Omega} \mA(\theta,:)^\top \mA(\theta,:) \,\dd\rho(\theta).
\]

Optimization is then conducted directly over the probability measure space, bypassing the limitations of finite discretization.

\medskip

Taken together, these directions highlight a broader landscape of experimental design beyond RNLA. In addition to embracing nonlinearity through adaptive and interactive strategies that tightly couple experimental design with reconstruction, and advancing toward infinite-dimensional formulations via tools from optimal transport, gradient flows, and mean-field theory, other approaches are also expected to play an important role. Vice versa,  experimental design methods have successfully been customized to several machine learning tasks  \cite{krause2008near, ash2021gone, cohn1996active}, prompting towards a more efficient use of data and training capacities. We anticipate that progress along these lines—at the interface of applied mathematics, machine learning, and computational science—will significantly broaden the scope and amplify the impact of experimental design methodologies in the coming decade.

\bibliographystyle{plain}
\bibliography{references.bib} 
\end{document}